\documentclass[12pt,reqno]{amsart}
\usepackage[pdfborder={0 0 0.5 [3 2]}, plainpages=false]{hyperref}%
\usepackage[left=1in,right=1in,top=1in,bottom=1in]{geometry}%
\usepackage[citation-order]{amsrefs}%
\usepackage{amsmath}
\usepackage{enumerate}
\usepackage{amssymb}                
\usepackage{amsfonts}
\usepackage{amsthm}
\usepackage{bbm}
\usepackage[table,xcdraw]{xcolor}
\usepackage{float}
\usepackage{mathtools}
\usepackage{cool}
\usepackage{graphicx,epsfig}

\usepackage[capitalize,nameinlink]{cleveref}
\crefname{section}{section}{sections}
\crefname{subsection}{subsection}{subsections}
\Crefname{section}{Section}{Sections}
\Crefname{subsection}{Subsection}{Subsections}

\Crefname{figure}{Figure}{Figures}

\crefformat{equation}{\textup{#2(#1)#3}}
\crefrangeformat{equation}{\textup{#3(#1)#4--#5(#2)#6}}
\crefmultiformat{equation}{\textup{#2(#1)#3}}{ and \textup{#2(#1)#3}}
{, \textup{#2(#1)#3}}{, and \textup{#2(#1)#3}}
\crefrangemultiformat{equation}{\textup{#3(#1)#4--#5(#2)#6}}%
{ and \textup{#3(#1)#4--#5(#2)#6}}{, \textup{#3(#1)#4--#5(#2)#6}}{, and \textup{#3(#1)#4--#5(#2)#6}}

\Crefformat{equation}{#2Equation~\textup{(#1)}#3}
\Crefrangeformat{equation}{Equations~\textup{#3(#1)#4--#5(#2)#6}}
\Crefmultiformat{equation}{Equations~\textup{#2(#1)#3}}{ and \textup{#2(#1)#3}}
{, \textup{#2(#1)#3}}{, and \textup{#2(#1)#3}}
\Crefrangemultiformat{equation}{Equations~\textup{#3(#1)#4--#5(#2)#6}}%
{ and \textup{#3(#1)#4--#5(#2)#6}}{, \textup{#3(#1)#4--#5(#2)#6}}{, and \textup{#3(#1)#4--#5(#2)#6}}

\crefdefaultlabelformat{#2\textup{#1}#3}

\def\R{{\mathbb R}}

\def\C{{\mathbb C}}
\def\Z{{\mathbb Z}}

\DeclareMathOperator{\spn}{span}
\DeclareMathOperator{\ran}{range}

\graphicspath{ {images/} }

\newtheorem{theorem}{Theorem}
\newtheorem{corollary}{Corollary}

\newtheorem{hypothesis}{Hypothesis}
\newtheorem{remark}{Remark}

\begin{document}

\title{Stationary multi-kinks in the discrete sine-Gordon equation}

\author{Ross Parker}
\address{Department of Mathematics, Southern Methodist University, 
Dallas, TX 75275, USA}
\email{rhparker@smu.edu}

\author{P.\,G. Kevrekidis} 
\address{Department of Mathematics and Statistics, University of Massachusetts, Amherst MA 01003, USA}
\email{kevrekid@math.umass.edu}

\author{Alejandro Aceves}
\address{Department of Mathematics, Southern Methodist University, 
Dallas, TX 75275, USA}
\email{aaceves@smu.edu}

\begin{abstract}
	We consider the existence and spectral stability of static multi-kink structures in the discrete sine-Gordon equation, as a representative example of the family of discrete Klein-Gordon models. The multi-kinks are constructed using Lin's method from an alternating sequence of well-separated kink and antikink solutions. We then locate the point spectrum associated with these multi-kink solutions by reducing the spectral problem to a matrix equation. For an $m$-structure multi-kink, there will be $m$ eigenvalues in the point spectrum near each eigenvalue of the primary kink, and, as long as the spectrum of the primary kink is imaginary, the spectrum of the multi-kink will be as well. We obtain analytic expressions for the eigenvalues of a multi-kink in terms of the eigenvalues and corresponding eigenfunctions of the primary kink, and these are in very good agreement with numerical results. We also perform numerical time-stepping experiments on perturbations of multi-kinks, and the outcomes of these simulations are interpreted using the spectral results.
\end{abstract}

\maketitle

\section{Introduction}

The $1+1$-dimensional nonlinear Klein-Gordon
models have been one of the preeminent platforms 
where ideas from integrable, as well as near-integrable, theory of nonlinear
partial differential equations and the corresponding 
solitary waves have been developed. This
is by now evidenced by numerous monographs~\cites{eilbeck,dauxois},
specialized books for these models~\cites{braun2004,SGbook,p4book} 
and reviews~\cites{kivsharmalomed,braun1998}.
Both continuum and discrete models of this kind
have been of interest to physical applications, as well
as to applied analysis, and their comparison as regards
the transition from discrete to continuum has been
of interest in its own right~\cite{SGchapter}.

The discrete Klein-Gordon equation
\begin{equation*}
\ddot{u}_n = d (\Delta_2 u)_n - f(u_n)
\end{equation*}
describes the dynamics of an infinitely long, one-dimensional lattice of particles which are harmonically coupled to their neighbors through the discrete second difference operator $\Delta_2$ and are subject to an external, nonlinear, onsite potential $P(u)$ such that $f(u) = P'(u)$ \cite{Karachalios}. The quantity $u_n$ represents the displacement of the particle at site $n$ in the lattice, $d$ is the strength of the nearest neighbor coupling, and the dot denotes the derivative with respect to the time $t$. 

A specific example is the discrete sine-Gordon equation
\begin{equation}\label{eq:dSG}
	\ddot{u}_n = d (\Delta_2 u)_n \pm \sin(u_n),
\end{equation}
in which the external potential is periodic. This equation is also known as the Frenkel-Kontorova model, and was introduced in 1938 to describe the dynamics of a crystal lattice near a dislocation core \cites{braun1998,braun2004}. This equation has since been used in numerous applications, including a mechanical model for a chain of pendula coupled with elastic springs \cites{Scott1969,english}, arrays of Josephson junctions \cites{Ustinov1992,Floria1998}, and DNA dynamics \cites{Yomosa1983,Yakushevich1998,DeLeo2011}. (See \cite{braun2004}*{Chapter 2} for more physical applications of this model). Another example of  substantial interest is the discrete $\phi^4$ model
\begin{equation}\label{eq:dphi4}
	\ddot{u}_n = d (\Delta_2 u)_n + 2 u_n(1 - u_n^2),
\end{equation}
which has a double well external potential, and has applications to conducting polymers \cite{heeger1988}. In addition, the
latter model has been a central point of focus
as concerns the dynamics of discrete 
breathers, i.e. time-periodic and exponentially localized solutions in space~\cites{tsironis,FLACH20081}.

\begin{figure}
	\begin{center}
	\begin{tabular}{ccc}
	\includegraphics[width=5cm]{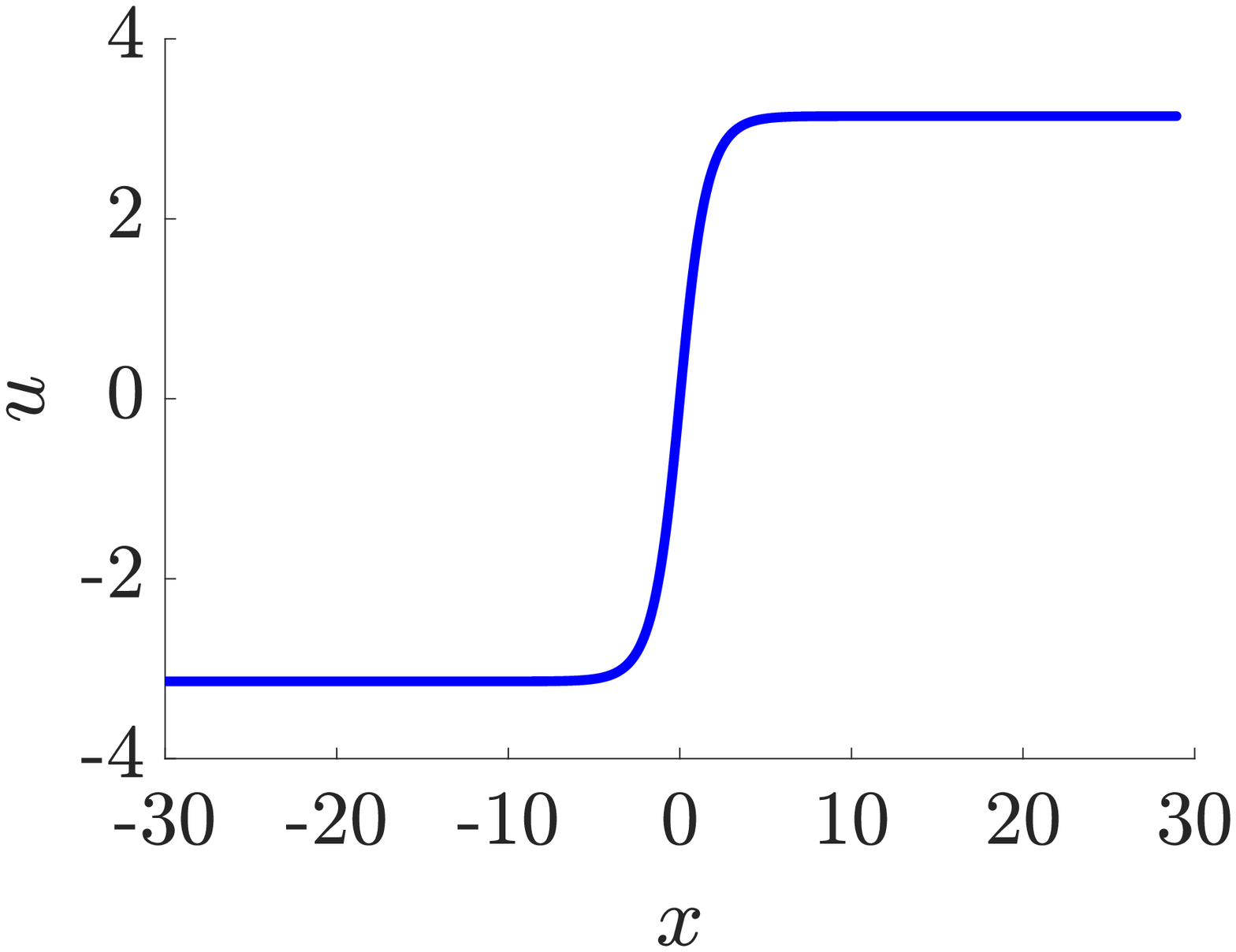}	&
	\includegraphics[width=5cm]{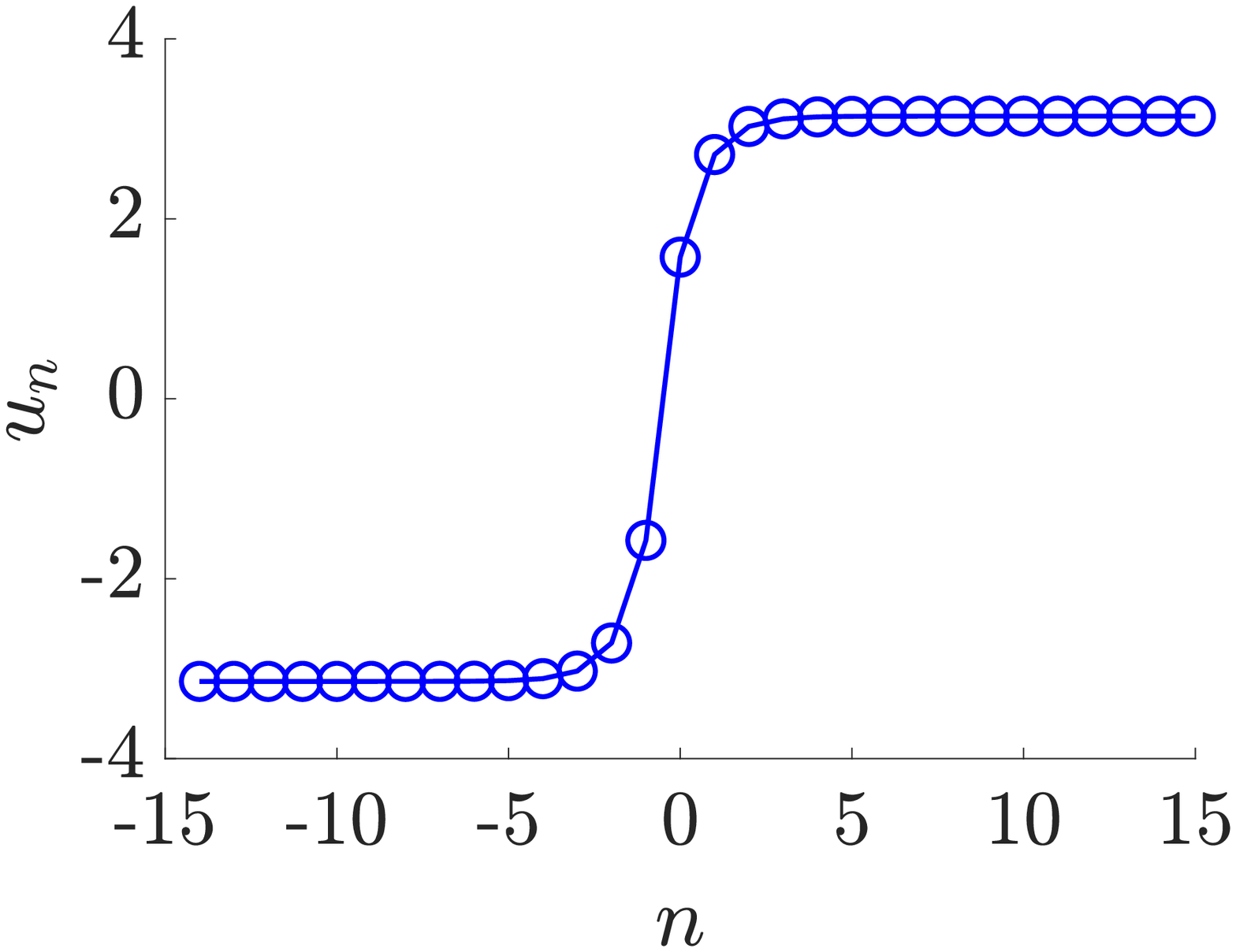} &
	\includegraphics[width=5cm]{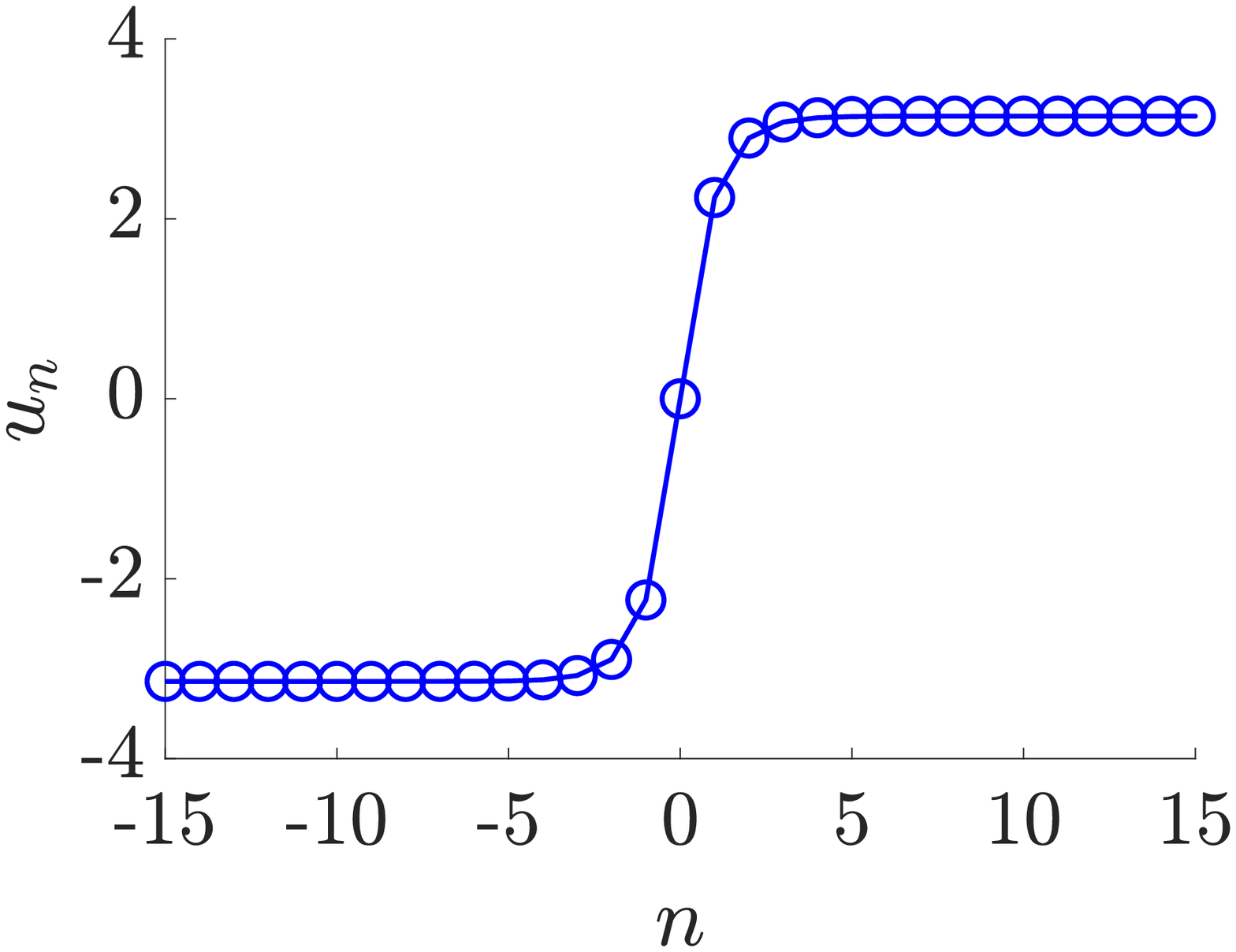}
	\end{tabular}
	\end{center}
	\caption{Kink solution to the sine-Gordon equation (left). Inter-site spectrally stable (center) and onsite 
	spectrally unstable (right) kinks of the discrete sine-Gordon equation for $d = 0.5$ and $n=60$ lattice points. }
	\label{fig:SGkinks}
\end{figure}

Equation \cref{eq:dSG} is the discrete analogue of the 
continuum sine-Gordon PDE
\begin{equation}\label{eq:SG}
	u_{tt} = u_{xx} \pm \sin(u),
\end{equation}
which has many physical (including, e.g.,
fluxons in Josephson junctions, charge density
waves in quasi-1d conducting materials, among many others)
and biological applications \cites{SGbook,Ivancevic2013} and has been extensively studied both in the mathematical and physics literature due to the fact that it is integrable via the inverse scattering transform \cite{SolitonBook1}. 
We note that if $u(x)$ is a solution to \cref{eq:SG} with the ``minus'' nonlinearity, $u(x) - \pi$ is a solution with the ``plus'' nonlinearity. The same holds for \cref{eq:dSG}. For symmetry reasons, which will be explained in \cref{sec:bg}, we will consider only the {``plus''} nonlinearity here.
Of particular interest are coherent structures such as kinks (\cref{fig:SGkinks}, left panel), exponentially localized stationary solutions which 
are heteroclinic orbits that connect two adjacent minima of the potential $P(u)$, and also continuum variants of the breather
solutions. Analytical waveforms are available via the 
integrable theory for both kinks and breathers in the continuum sine-Gordon equation (see, for example, \cite{SGchapter}). Kink  solutions exist for the continuum $\phi^4$ model as well, although that equation is not integrable \cites{SGbook,KevrekidisWeinstein2000}
and the same is generically true for Klein-Gordon
models
with multiple degenerate energy minima.
Regular breather waveforms do not exist
in the continuum $\phi^4$ model (and generically in continuum
Klein-Gordon settings aside from the integrable
sine-Gordon case) as is known from the work of~\cite{segur}. 

When we move from the continuum realm to the discrete, there are two distinct kink solutions: intersite kinks (\cref{fig:SGkinks}, which connect two adjacent minima of the potential $P(u)$ directly, and onsite kinks (\cref{fig:SGkinks}, right panel), which connect these states via a local maximum of $P(u)$ which lies between the two minima~\cite{peyrard}. 
This is similar to the discrete nonlinear Schr{\"o}dinger equation (DNLS), where there is both a site-centered and an intersite-centered (pulse) soliton solution~\cite{Kevrekidis2009}. Corresponding to each kink solution is an antikink solution which connects the two equilbria in the opposite order. By reversibility, if $k(n)$ is a kink solution, then $k(-n)$ is an antikink. We will be concerned only with static kink solutions here, although we note that there has been much interest in both moving kinks \cites{Aigner2003,Iooss2006,Cisneros2008} and breather solutions (see, for example, \cites{SGbook} as well as \cites{Pelinovsky2012,Cuevas2011} for results on multi-site breathers) in discrete Klein-Gordon lattices.

In this paper, we look at multi-kinks, which are coherent structures resembling a sequence of alternating kinks and antikinks spliced together end-to-end. 
The study of complex coherent structures formed by joining multiple copies of a simpler coherent structure has a rich mathematical history (see \cite{Sandstede1998}, and the references therein). In \cite{Sandstede1998}, for example, the existence and stability of multi-pulse solutions to semilinear parabolic equations, such as reaction diffusion equations, is determined using Lin's method \cites{Lin1990,Lin2008}, an implementation of the Lyapunov-Schmidt reduction. This method constructs multi-pulses by splicing together multiple, well-separated copies of a single pulse using small remainder functions, and it is also used to reduce the PDE eigenvalue problem to a matrix eigenvalue problem. These techniques have been recently extended to Hamiltonian lattice systems, including DNLS \cite{Parker2020}. In the case of DNLS, multi-pulse solutions exist on the lattice which do not exist in the continuum equation (although multi-pulses do exist in higher order NLS models \cite{Parker2021}). 

For the discrete Klein-Gordon equation, we follow a similar approach to \cite{Parker2020} and use a discrete adaptation of Lin's method \cite{Knobloch2000} to construct multi-kinks from a sequence of well-separated kinks and antikinks. As with multi-pulses in DNLS, these multi-kinks do not exist in the continuum equation. 
Indeed, it is the so-called Peierls-Nabarro barrier~\cite{peyrard}, the
local effective potential due to discreteness, that makes
such configurations possible. Otherwise, similarly to the
DNLS case~\cite{Kapitula2001a}, the structures would purely 
interact exponentially through their tails, being unable
to form stationary patterns involving multiple waves.
We also use Lin's method to reduce the spectral problem to a matrix equation. The most notable difference here is that the system on the lattice is no longer translation invariant, thus there is not an eigenvalue at 0. Instead of locating eigenvalues in a small ball around the origin, this method locates eigenvalues of the multi-kink in the neighborhood of each eigenvalue of the primary kink. In particular, this means that there will be a different matrix reduction corresponding to each eigenvalue of the primary kink; computing each of these matrices requires knowing (or computing numerically) the point spectrum and the corresponding eigenfunctions of the primary kink. We show that for an $m$-component multi-kink, there are $m$ eigenvalues in the point spectrum near each eigenvalue of the primary kink. We obtain analytic expressions for these eigenvalues in terms of the tails of the eigenfunctions of the primary kink, in a similar fashion to the expressions for DNLS which depend on the tails of the primary pulse \cite{Parker2020}. These expressions are in good agreement with numerical computations for intermediate values of the coupling parameter $d$. Finally, we perform timestepping experiments on perturbations of the primary kink and multi-kinks to illustrate the role of the point spectrum and its corresponding eigenfunctions in explaining the evolution of these perturbations. We believe that this offers a 
systematic understanding of such multiwave patterns, based
on the existence and stability properties of the corresponding
building block, namely the single kink (or antikink).

This paper is organized as follows. In \cref{sec:bg}, we present the mathematical background for the discrete Klein-Gordon equation, together with a reformulation of the existence and eigenvalue problems using a spatial dynamics approach. The main results concerning the existence and spectrum of multi-kinks are then given in \cref{sec:multikink}; the proofs of these results are deferred to the end of the paper. In \cref{sec:numerics}, we present numerical results which corroborate the main theorems, as well as results of timestepping simulations. The proofs of \cref{th:KaKexists} and \cref{th:stability} are given as appendices following a brief concluding section.

\section{Mathematical background}\label{sec:bg}

We will consider the discrete Klein-Gordon equation with onsite nonlinearity $f(u)$
\begin{equation}\label{eq:KG}
\ddot{u}_n = d (\Delta_2 u)_n - f(u_n),
\end{equation}
where $(\Delta_2 u)_n = u_{n+1} - 2 u_n + u_{n-1}$ is the discrete second difference operator, and $f(u) = P'(u)$ for a smooth potential function $P(u)$. The nonlinearity $f$ has the following three properties:
\begin{enumerate}[(i)]
	\item $f(u)$ is an odd function with $f'(0) < 0$. This implies that $f(0) = 0$, so 0 is an equilibrium of \cref{eq:KG}.
	\item There is a pair of nonzero equilibria $\pm u^*$ with $f(\pm u^*) = 0$ and $f'(u^*) = f'(-u^*) > 0$.
	\item There are no other equilibria in $[-u^*, u^*]$.
\end{enumerate}
Important versions include the discrete sine-Gordon equation, where $f(u) = -\sin(u)$ and $P(u) = 1 + \cos(u)$, and the $\phi^4$ model, where $f(u) = -u(1-u^2)$ and $P(u) = \frac{1}{4}(1-u^2)^2$.
We note that the discrete sine-Gordon equation is typically written with the nonlinearity $f(u) = \sin u$, in which case the pair of stable
equilibria in (ii) are at 0 and $2 \pi$, and the unstable one between them
is at $\pi$. This is identical to the situation considered here except for a shift of $u_n$ by $-\pi$. Choosing $f$ to have odd symmetry greatly simplifies the analysis, and allows the result to apply to both the discrete sine-Gordon equation and the $\phi^4$ model.
Equation \cref{eq:KG} is Hamiltonian, with energy given by \cite{KevrekidisWeinstein2000}
\begin{equation}\label{eq:H}
	\mathcal{H}(u) = \sum_{n=-\infty}^\infty 
	\left( \frac{1}{2} (\dot{u}_n)^2 + \frac{d}{2} (u_{n+1} - u_n)^2 + P(u_n) \right).
\end{equation}

The equilibrium solutions we will study are standing waves which satisfy 
\begin{equation}\label{eq:KGeq}
d (\Delta_2 u)_n - f(u_n) = 0.
\end{equation}
Linearization about an equilibrium solution $u_n$ yields the eigenvalue problem
\begin{equation}\label{eq:KGevp1}
d (\Delta_2 v)_n - f'(u_n)v_n = \lambda^2 v_n.
\end{equation}
Letting $\omega = \lambda^2$, we obtain the eigenvalue problem for $\omega$
\begin{equation}\label{eq:evp}
d (\Delta_2 v)_n - f'(u_n)v_n = \omega v_n.
\end{equation}
The eigenvalues are given by $\lambda = \pm \sqrt{\omega}$. Equation \cref{eq:evp} has the form of an infinite dimensional matrix problem \cite{Balmforth2000}. Since that matrix is real and symmetric, the linear operator defined by the LHS of \cref{eq:evp} is self-adjoint, thus $\omega$ must be real. This implies that the eigenvalues $\lambda$ must be either real or purely imaginary pairs. In particular, spectral instabilities can only develop when a pair of eigenvalues passes through the origin \cite{Balmforth2000}.

Using a spatial dynamics approach as in \cite{Parker2020}, let $u_n$ be an equilibrium solution to \cref{eq:KGeq}, and let $U(n) = (u(n), \tilde{u}(n)) = (u_n, u_{n-1})$. Then equation \cref{eq:KGeq} is equivalent to the lattice dynamical system in $\R^2$
\begin{equation}\label{eq:dynEq}
U(n+1) = F(U(n)),
\end{equation}
where
\[
F\begin{pmatrix}u \\ \tilde{u} \end{pmatrix} =
\begin{pmatrix}2u - \tilde{u} + \frac{1}{d}f(u) \\
u
\end{pmatrix}.
\]
Equation \cref{eq:dynEq} has three fixed points of interest at 0 and $S^\pm = (\pm u^*, \pm u^*)^T$. (It may in fact have more, depending on the specific form of $f(u)$; this is the case for the discrete sine-Gordon equation). Linearizing about the fixed points $S^\pm$, we obtain the matrix 
\[
D F(S^\pm)=
\begin{pmatrix}2 + \frac{1}{d}f'(\pm u^*) & -1 \\ 1 & 0
\end{pmatrix}.
\]
Since $f'(\pm u^*) > 0$, this has a pair of eigenvalues $\{ r, 1/r \}$ with $r > 0$, where
\begin{equation}\label{eq:r}
r = \frac{1}{2d}\left( f'(u^*) + 2d + \sqrt{f'(u^*)(f'(u^*) + 4d)} \right).
\end{equation}
Thus $S^\pm$ are hyperbolic saddle equilibria of the lattice dynamical system \cref{eq:dynEq}. The origin is a nonhyperbolic equilibrium which has a pair of eigenvalues on the unit circle in the complex plane. We take the existence of a stable, symmetric kink (stationary front) as a hypothesis. From the spatial dynamics perspective, this kink solution is a heteroclinic orbit connecting the saddle at $S^-$ to the saddle at $S^+$.

\begin{hypothesis}\label{hyp:kinkexists}
There exists a kink solution $K(n) = (k(n),\tilde{k}(n))$ to \cref{eq:dynEq} which connects the unstable manifold $W^u(-u^*, -u^*)$ and the stable manifold $W^s(u^*, u^*)$. These manifolds intersect transversely in $\R^2$. 
The kink has the odd symmetry $k(-n) = -k(n-1)$. Finally, the kink $K(n)$ is a minimizer of the $t$-independent energy functional
\begin{equation}\label{eq:energyfunctional}
h[u] = \sum_{n=-\infty}^\infty 
\left( \frac{d}{2} (u_{n+1} - u_n)^2 + P(u_n) \right)
\end{equation}
{among the class of heteroclinic connections between $S^-$ and $S^+$.}
\end{hypothesis}

\begin{remark}
The kink $k(n)$ from \cref{hyp:kinkexists} has the odd symmetry $k(-n) = -k(n-1)$, and is known as an intersite kink, since it does not involve the equilibrium at 0. For the sine-Gordon equation, an intersite kink at the anti-continuum (AC) limit 
{$d = 0$} is given by $(\dots, -\pi, -\pi, \pi, \pi, \dots)$. By contrast, an onsite kink, which has the odd symmetry $k(-n) = -k(n)$, involves the equilibrium at 0. For the sine-Gordon equation, an onsite kink at the AC limit is $(\dots, -\pi, -\pi, 0, \pi, \pi, \dots)$. 
\end{remark}

Since $K(n)$ is a minimizer of the energy functional \cref{eq:energyfunctional}, the spectrum of $K(n)$ lies on the imaginary axis~\cite{KevrekidisWeinstein2000}*{Section 2.1.6}.
For specific nonlinearities $f(u)$, including those associated with the discrete sine-Gordon equation and the $\phi^4$ model, the existence of a symmetric, intersite kink which is a minimizer of \cref{eq:energyfunctional} is known (see \cites{KevrekidisWeinstein2000,SGchapter} and references therein). On the other
hand, the onsite kink will in general be unstable. For the sine-Gordon equation, for example, the onsite kink has a pair of real eigenvalues $\pm \lambda$ and is thus unstable \cite{Kapitula2001}*{Theorem 4.4}. In addition, the onsite kink can be shown to be unstable for the sine-Gordon equation when $d < 1/4$ and the  $\phi^4$ model when $d < 1/2$ using Gerschgorin’s theorem \cite{SGchapter}. 
While the former result is an asymptotic one, valid
in the vicinity of the continuum limit, the latter is a
rigorous one, but {\it only} valid for the above-mentioned
interval of $d$ in the vicinity of the AC limit
of $d=0$. Hence, the two results are complementary to each other.
Since $f(u)$ is an odd function, if $u_n$ is a solution to \cref{eq:KGeq}, so is $-u_n$. Thus for every kink solution $K(n)$ to \cref{eq:dynEq} there is a corresponding antikink solution $\tilde{K}(n) = -K(n)$.

The spectrum of the primary kink solution $K(n) = (k(n),\tilde{k}(n))$ can be decomposed into two disjoint sets: the point spectrum consists of isolated eigenvalues for which the corresponding eigenfunction is in $\ell^2(\Z)$, and the continuous spectrum which consists of bounded, oscillatory modes. Following \cite{KevrekidisWeinstein2000}, the continuous spectrum depends only on the background state of the system and consists of the two symmetric intervals on the imaginary axis
\begin{equation}\label{eq:contspec}
	\sigma_{\text{cont}} = \pm i \left[\sqrt{f'(u^*)}, \sqrt{f'(u^*) + 4d}\right].
\end{equation}
In particular, there is a gap in the continuous spectrum $i\left(-\sqrt{f'(u^*)},\sqrt{f'(u^*)}\right)$ which contains the origin. The continuous spectrum will be the same for the linearization about any equilibrium solution involving the asymptotics of $u \rightarrow \pm u^*$.

The continuum Klein-Gordon equation $u_{tt} = u_{xx} - f(u)$ has an eigenvalue at 0 due to translation invariance which is referred to as the Goldstone mode. Since the discrete Klein-Gordon equation does not possess any continuous symmetries, there will be no eigenvalues at the origin. Instead, there will be a symmetric pair of eigenvalues, which is either real or purely imaginary. For the intersite kink we are considering, this pair will be imaginary since the entire spectrum is imaginary (it will be real for the onsite kink). This pair of eigenvalues is often termed Goldstone modes by extension of the Goldstone mode of the continuum equation, since these eigenvalues approach the origin as the discrete equation approaches the continuum limit, i.e. $d \rightarrow \infty$. For specific nonlinearities $f(u)$ and certain values of the coupling parameter $d$, there may be additional eigenvalues, known as internal modes, which lie between the bands of the continuous spectrum. For the kinks which are minimizers of \cref{eq:energyfunctional}, these will also be purely imaginary. (See \cites{cretegny,KevrekidisWeinstein2000} for a discussion of the spectrum of the kink solution, including the internal modes, for the discrete sine-Gordon equation and the $\phi^4$ model). We make the additional hypothesis that all point spectrum for the primary kink $K(n)$ lies in the continuous spectrum gap.

{
\begin{hypothesis}\label{hyp:pointspecgap}
$|\lambda| < \sqrt{f'(u^*)}$ for all eigenvalues (point spectrum) $\lambda$ of the primary kink $K(n)$.
\end{hypothesis}
}

\section{Multi-kinks}\label{sec:multikink}

We can construct multi-kink solutions by joining together an alternating sequence of kinks and antikinks in an end-to-end fashion. By the symmetry of $f(u)$, $u_n$ is an equilibrium solution if and only if $-u_n$ is, thus we can always without loss of generality begin with a kink solution. We will characterize a multi-kink in the following way. Let $m > 1$ be the total number of kinks and antikinks. Let $N_i$ ($i = 1, \dots, m-1$) be the distances (in lattice points) between consecutive kinks/antikinks.
(As a mark of the position of the center of the kinks/antikinks,
we use the position of the corresponding zero-crossing).
We seek a solution which can be written piecewise in the form 
\begin{equation}\label{eq:Upiecewise}
\begin{aligned}
U_i^-(n) &= c_i K(n) + \tilde{U}_i^-(n) && n \in [-N_{i-1}^-, 0] && \quad i = 1, \dots, m\\
U_i^+(n) &= c_i K(n) + \tilde{U}_i^+(n) && n \in [0, N_i^+] && \quad i = 1, \dots, m,
\end{aligned}
\end{equation}
where $c_i = (-1)^{i+1}$, $N_i^+ = \lfloor \frac{N_i}{2} \rfloor$, $N_i^- = N_i - N_i^+$, and $N_0^- = N_m^+ = \infty$. 
We also define
\begin{equation}\label{defN}
N = \frac{1}{2} \min\{ N_i \}
\end{equation}
as a characteristic distance, which will be used in the estimates of the remainder terms in \cref{eq:Uestimates}.
The individual pieces $U_i^\pm(n)$ are joined together end-to-end as in \cites{Sandstede1998,Knobloch2000,Parker2020} to create the multi-kink solution $U(n)$, which can be written in piecewise form as
\begin{equation}
\begin{aligned}
U(n) &= \begin{cases}
U_i^-\left( n - \sum_{j=1}^{i-1}N_j \right) & \sum_{j=1}^{i-1}N_j - N_{i-1}^- + 1 \leq n \leq \sum_{j=1}^{i-1}N_j \\
U_i^+\left( n - \sum_{j=1}^{i-1}N_j \right) & \sum_{j=1}^{i-1}N_j + 1 \leq n \leq \sum_{j=1}^{i-1}N_j + N_i^+
\end{cases}
&& i = 1, \dots, m,
\end{aligned}
\end{equation}
where we define $\sum_{j=1}^0 N_j = 0$. Since we are taking $N_0^- = N_m^+ = \infty$, this formula makes sense for the two end pieces $U_1^-(n)$ and $U_m^+(n)$.
The functions $\tilde{K}_i^\pm(n)$ in \cref{eq:Upiecewise} are remainder terms, which will be small. We then have the following existence theorem. The proof is deferred until \cref{sec:proof1}.

\begin{theorem}\label{th:KaKexists}
Assume \cref{hyp:kinkexists}. Then there exists a positive integer $N_0$ with the following property. For all $m > 1$ and distances $N_i \geq N_0$, there exists a unique solution $U(n)$ which is composed of $m$ alternating kinks and antikinks and can be written piecewise in the form \cref{eq:Upiecewise}. For the remainder terms $\tilde{U}_i^\pm(n)$, we have the estimates
\begin{equation}\label{eq:Uestimates}
\begin{aligned}
\|\tilde{U}_i^\pm\| &\leq C r^{-N} \\
| \tilde{U}_i^-(n) | &\leq C r^{-N_{i-1}^-} r^{-(N_{i-1}^- + n)} && \qquad n = 2, \dots, m\\
|\tilde{U}_i^+(n)| &\leq C r^{-N_i^+} r^{-(N_i^+ - n)} && \qquad n = 1, \dots, m-1 \\
| \tilde{U}_1^-(n) | &\leq C r^{-2N} r^{n} \\
|\tilde{U}_m^+(n)| &\leq C r^{-2N} r^{-n} .
\end{aligned}
\end{equation}
\end{theorem}

\begin{remark}\label{rem:SGmulitkinks}
For the sine-Gordon equation, equation \cref{eq:dynEq} has saddle equilibria at $n \pi$ for all odd integers $n$. More general multi-kinks can be constructed comprising kinks and antikinks which link any adjacent saddle equilibria. A specific example is a double kink, where a kink  connecting $-\pi$ to $\pi$ is followed by a kink connecting $\pi$ to $3 \pi$. The existence of these general multi-kinks is a straightforward adaptation of \cref{th:KaKexists}.
\end{remark}

To determine the eigenvalues of the linearization about a multi-kink, we again take a spatial dynamics approach. We rewrite \cref{eq:evp} as a lattice dynamical system by taking $V(n) = (v(n), \tilde{v}(n)) = (v_n, v_{n-1})$. Then \cref{eq:evp} is equivalent to the lattice dynamical system in $\R^2$
\begin{equation}\label{eq:EVPdyneq}
V(n+1) = D F( U(n) )V(n) + \omega B V(n),
\end{equation}
where
\[
B = \frac{1}{d}
\begin{pmatrix}1 & 0 \\ 0 & 0
\end{pmatrix}.
\] 
For the primary kink $K(n) = (k_n, \tilde{k}_n)$, let $v_0(n)$ be an eigenfunction with corresponding eigenvalue $\lambda_0$, let $\omega_0 = \lambda_0^2$, and let $V_0(n) = (v_0(n), \tilde{v}_0(n)) = (v_0(n), v_0(n-1))$. Then $V_0(n)$ solves the equation
\begin{equation}\label{eq:V0eq}
V_0(n+1) = D F( K(n) )V_0(n) + \omega_0 B V_0(n),
\end{equation}
which we rewrite as
\begin{equation}\label{eq:V0Aeq}
	V_0(n+1) = A(n; \omega_0) V_0(n),
\end{equation}
where
\begin{equation}\label{eq:Aomegaeq}
	A(n; \omega) = D F( K(n) )V_0(n) + \omega B.
\end{equation}
By the stable manifold theorem, 
\begin{equation}\label{eq:A0decay}
	|A(n; \omega_0) - A_0| \leq C r^{-|n|},
\end{equation}
where $A_0$ is the constant matrix
\begin{equation}
	A_0 = DF(S^+) + \omega_0 B.
\end{equation}
$A_0$ has eigenvalues $\{ r_0, 1/r_0 \}$, where
\begin{equation}\label{eq:r0}
r_0 = \frac{1}{2d}\left( f'(u^*) + \omega_0 + 2d + \sqrt{(f'(u^*)+ \omega_0)(f'(u^*) + \omega_0 + 4d)} \right).
\end{equation}
Since $\lambda_0$ is on the imaginary axis by \cref{hyp:kinkexists}, $\omega_0 < 0$, thus it follows from \cref{hyp:pointspecgap} that $A_0$ is hyperbolic.
It also follows from \cref{hyp:pointspecgap} that $1 < r_0 < r$, thus the eigenfunctions decay to 0 slower than the kink solution decays to the equilibria at $\pm S$.

For a multi-kink composed of $m$ components, each eigenvalue of the primary kink $K(n)$ will split into $m$ eigenvalues. The following theorem locates these eigenvalues for the multi-kink $U(n)$. The proof is deferred until \cref{sec:proof2}.

\begin{theorem}\label{th:stability}
Assume \cref{hyp:kinkexists} and \cref{hyp:pointspecgap}.
Let $U(n)$ be an $m-$component multi-kink constructed as in \cref{th:KaKexists} with distances $N_i$, and let $N = \frac{1}{2} \min\{ N_1, \dots, N_{m-1}\}$. Let $\pm \lambda_0$ be a pair of purely imaginary eigenvalues for the primary kink with corresponding eigenfunction $V_0(n)$, and let $\omega_0 = \lambda_0^2$. Then for $N$ sufficiently large, the multi-kink $U(n)$ has $m$ pairs of imaginary eigenvalues $\{\pm \lambda_0^1, \dots, \pm \lambda_0^m \}$ which are close to $\pm \lambda_0$ and are given by
\begin{align}\label{eq:lambdaj}
	\lambda_0^j = \sqrt{\omega_j}, \qquad
	\omega_j = \lambda_0^2 + \frac{d \mu_j}{M} + \mathcal{O}(r_0^{-3N}) && j = 1, \dots, m,
\end{align}
where $r_0$ is defined in \cref{eq:r0}, $\{ \mu_1, \dots, \mu_m \}$ are the real, distinct eigenvalues of the symmetric, tridiagonal matrix
\begin{equation}\label{eq:matrixA}
	A = \begin{pmatrix}
	0 & a_1 & & & \\
	a_1 & 0 & a_2 \\
	& a_2 & 0 & a_3 \\
	& & \ddots & \ddots & \\
	& & & & a_{m-1}  \\
	& & & a_{m-1} & 0  \\
	\end{pmatrix},
\end{equation}
with 
\begin{equation}\label{eq:ai}
	a_i = v_0(N_i^+)v_0(-N_i^- - 1) - v_0(N_i^+ - 1)v_0(-N_i^-),
\end{equation}
and $M$ is the Melnikov sum
\begin{equation}\label{eq:Minth}
	M = \sum_{n = -\infty}^{\infty} v_0(n)^2 = \| v_0 \|_{\ell^2}.
	\end{equation}
\end{theorem}

\begin{remark}\label{remark:unstablekink}
	The results of \cref{th:stability} hold as well for multi-kinks constructed using the unstable, onsite kink. In that case, the multi-kink $U(n)$ would have $m$ pairs of real eigenvalues close to each real eigenvalue pair $\pm \lambda_0$ of the onsite kink.
\end{remark}

\begin{remark}\label{remark:kk}
The estimates in \cref{th:stability}, and in particular the matrix \cref{eq:matrixA}, require knowledge of both the eigenvalue $\lambda_0$ of the primary kink and its corresponding eigenfunction $v_0(n)$. The matrix $A$ will be different for each eigenvalue $\lambda_0$ of the primary kink. In particular, note that the decay rate of the remainder term in \cref{eq:lambdaj} depends on $\lambda_0$ via $r_0$.
\end{remark}

\noindent For $m = 2$ and $m = 3$, we can compute the eigenvalues of $A$ exactly.

\begin{corollary}\label{corr:m23}
For $m = 2$, 
\begin{align*}
	\omega_1 &= \lambda_0^2 + \frac{d}{M}a_1 + \mathcal{O}(r_0^{-3N}) \\
	\omega_2 &= \lambda_0^2 - \frac{d}{M}a_1 + \mathcal{O}(r_0^{-3N}).
\end{align*}
For $m = 3$,
\begin{align*}
	\omega_1 &= \lambda_0^2 + \mathcal{O}(r_0^{-3N}) \\
	\omega_2 &= \lambda_0^2 + \frac{d}{M}\sqrt{a_1^2 + a_2^2} + \mathcal{O}(r_0^{-3N}) \\
	\omega_3 &= \lambda_0^2 - \frac{d}{M}\sqrt{a_1^2 + a_2^2} + \mathcal{O}(r_0^{-3N}).
\end{align*}
\end{corollary}

\begin{remark}
\cref{th:stability} can be easily adapted to the case of general multi-kinks in the discrete sine-Gordon equation (see \cref{rem:SGmulitkinks}). The matrix \cref{eq:matrixA} will have a similar form. The spectrum of all general multi-kinks will be purely imaginary as long as the multi-kink comprises only stable, intersite kinks.
\end{remark}

\section{Numerical results}\label{sec:numerics}

\begin{figure}
\begin{center}
\begin{tabular}{cc}
\includegraphics[width=7.5cm]{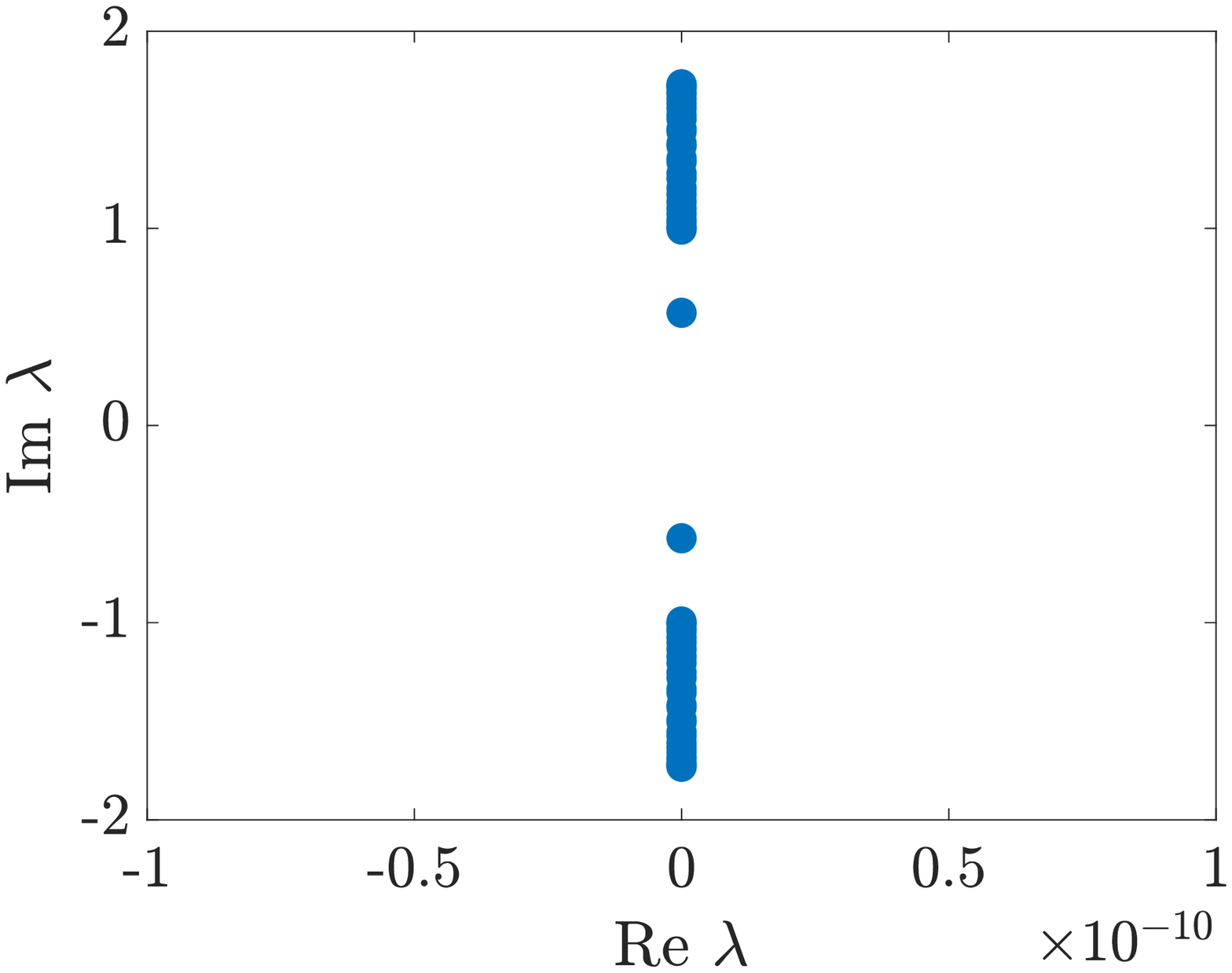}	&
\includegraphics[width=7.5cm]{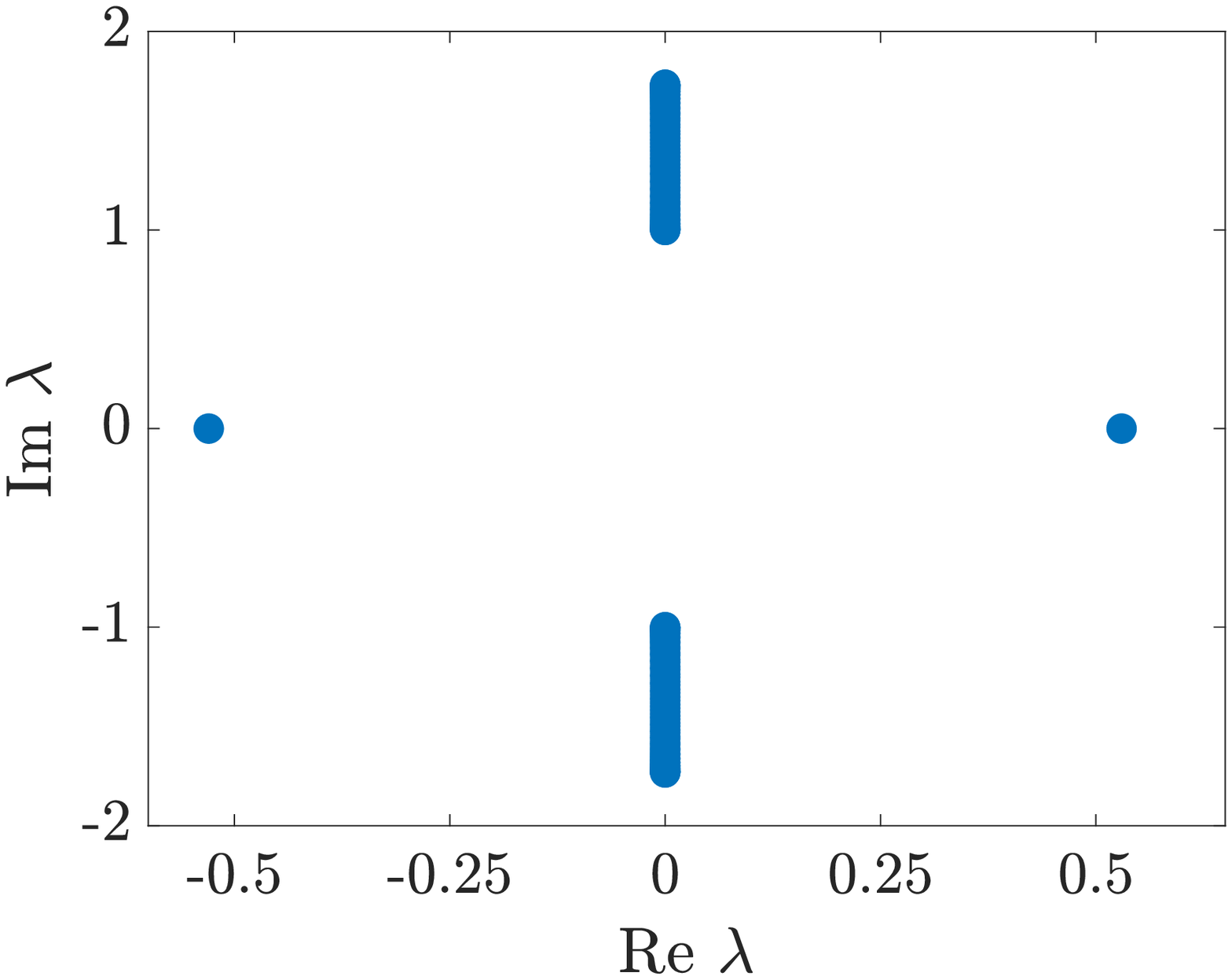}
\end{tabular}
\end{center}
\caption{Spectrum of the primary spectrally stable 
intersite kink (right) and of the spectrally unstable 
onsite kink (left) for the discrete sine-Gordon equation with $d = 0.50$.}
\label{fig:kinkspec}
\end{figure}

The results we present here are from the discrete sine-Gordon equation
$
\ddot{u}_n = d (\Delta_2 u)_n + \sin(u_n)
$.
We start by constructing the primary, intersite kink solution $k_n$ by using MATLAB for numerical parameter continuation from the anti-continuum (AC) limit ($d = 0$) in the coupling parameter $d$, starting with the solution $(\dots, -\pi, -\pi, \pi, \pi, \dots)$ (\cref{fig:SGkinks}, center). (The right panel of \cref{fig:SGkinks} shows the onsite kink). We then compute the spectrum of the linearization about the primary kink $k_n$ using MATLAB's \texttt{eig} function (\cref{fig:kinkspec}). As expected~\cites{Balmforth2000,KevrekidisWeinstein2000}, the spectrum of the intersite kink is imaginary, whereas the spectrum of the onsite kink contains a symmetric pair of real eigenvalues, thus it is unstable. Most of the remaining numerical results will only concern stable, intersite kinks.

\begin{figure}[H]
	\begin{center}
	\begin{tabular}{ccc}
	\includegraphics[width=5cm]{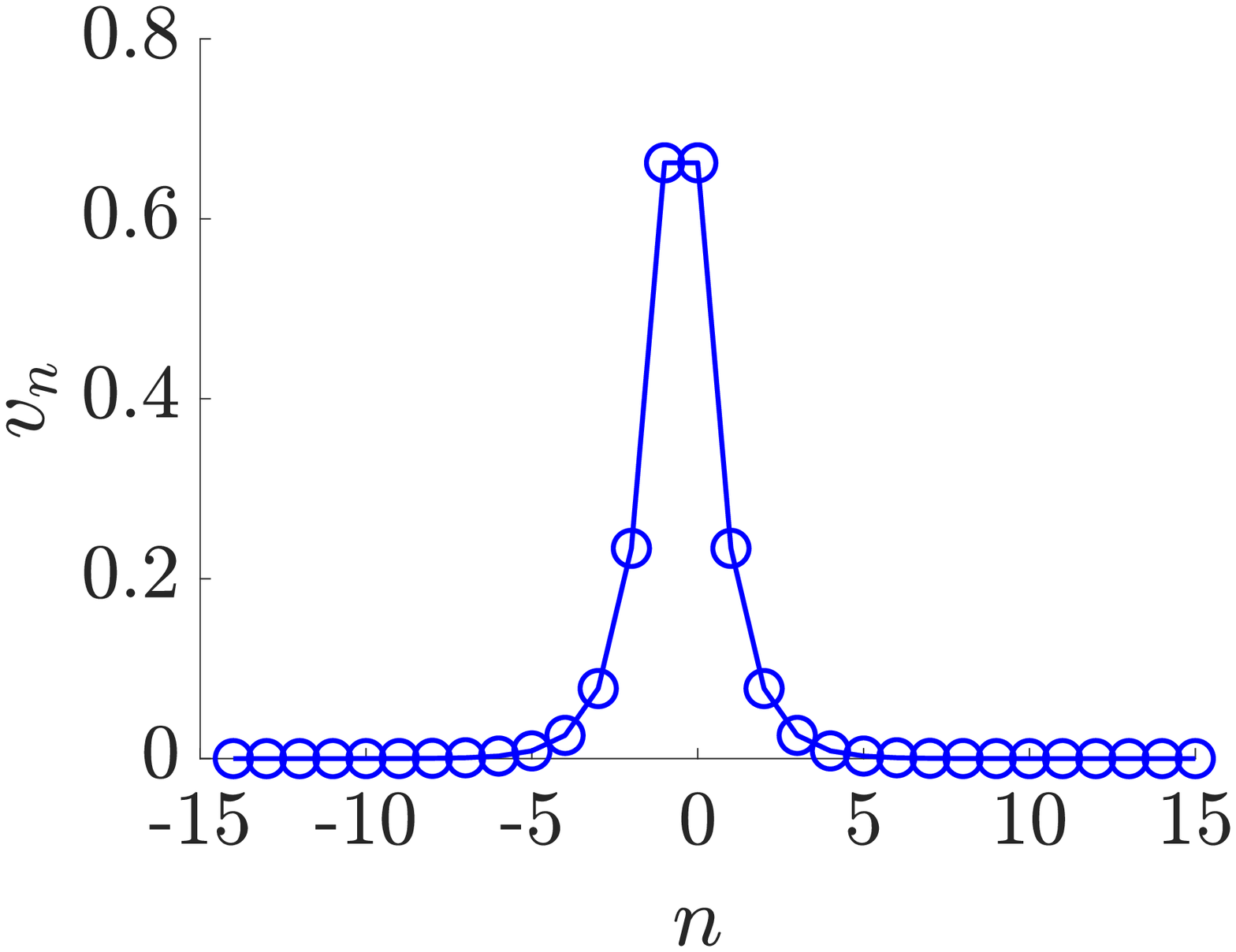} &
	\includegraphics[width=5cm]{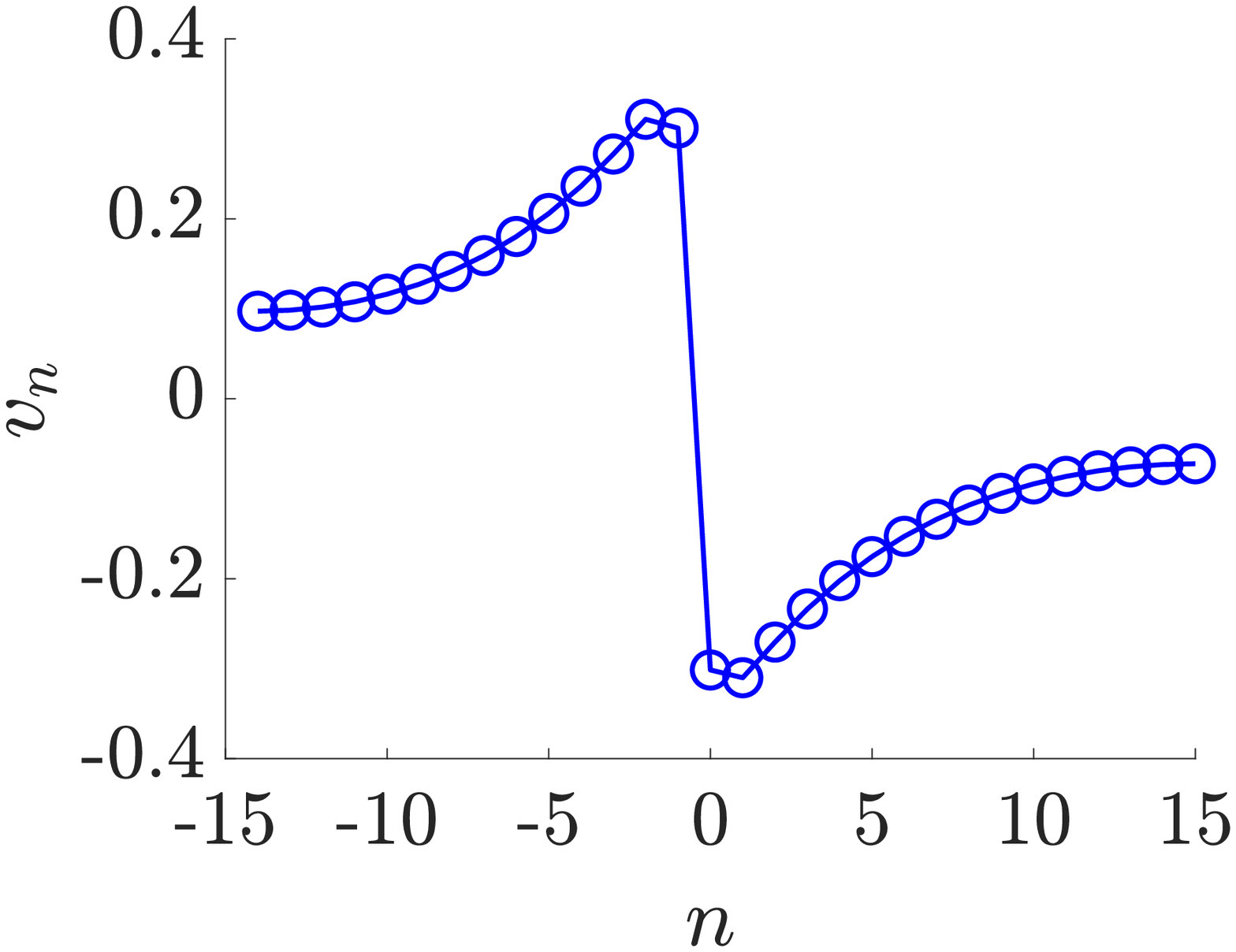} &
	\includegraphics[width=5cm]{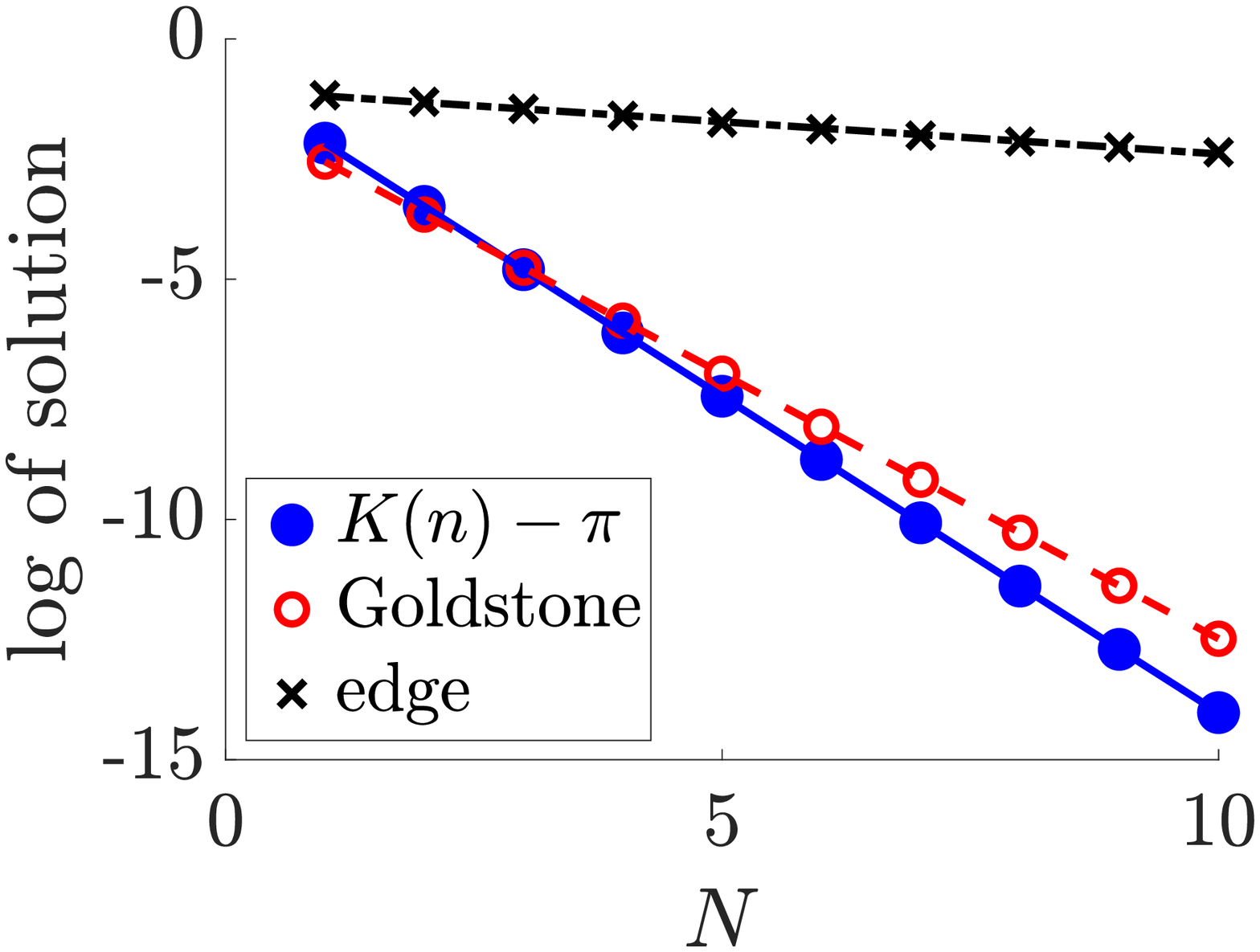}
	\end{tabular}
	\end{center}
	\caption{Goldstone mode eigenfunction ($\lambda = \pm 0.5718 i$, left) and edge mode eigenfunction ($\lambda = \pm 0.9941 i$, center). Semilog plot of decay of primary kink to $\pm \pi$ and decay of Goldstone and edge modes to 0 (right). Lines are least-squares linear regressions. $d = 0.50$ and $n =60$ lattice points. }
	\label{fig:kinkeig}
\end{figure}

The continuous spectrum lies within the interval given by \cref{eq:contspec}, and the pair of Goldstone mode eigenvalues is clearly visible in the continuous spectrum gap. For approximately $d > 0.265$, there is an additional internal mode eigenvalue for the intersite kink (not discernible in the left panel of \cref{fig:kinkspec}), which is known as an edge mode since it arises from the continuous spectrum (see \cite{KevrekidisWeinstein2000}*{Section 2.2}, noting that we are using $d$ in place of $d^2$ in that paper.) The eigenfunctions corresponding to the Goldstone mode and the edge mode are shown on the left and middle panels of of \cref{fig:kinkeig}. The semilog plot on the right panel shows the exponential decay of the primary kink to $\pi$ and the eigenfunctions to 0 as the lattice site index $n$ increases. These decay rates are predicted to be $r^{-|n|}$ for the primary kink, and $r_0^{-|n|}$ for the eigenfunctions, where $r$ is given by \cref{eq:r} and $r_0$ (which depends on the eigenvalue) is given by \cref{eq:r0}. The decay rates computed from the least squares linear regression lines have a relative error of order $10^{-4}$ for the kink and the Goldstone mode and a relative error of $10^{-2}$ for the edge mode.

\begin{figure}
	\begin{center}
	\includegraphics[width=16.5cm]{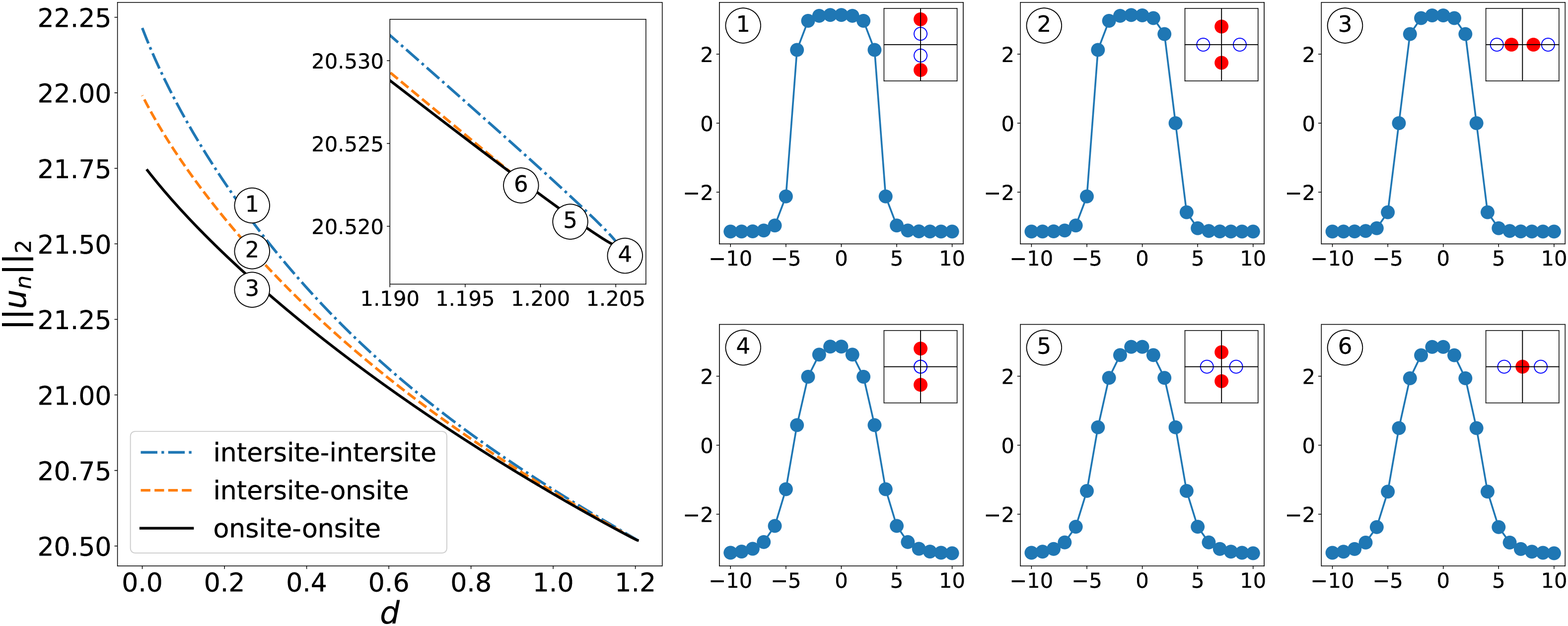}
	\end{center}
	\caption{The left panel shows the bifurcation diagram for a kink-antikink waveform
	with $N_1 = 8$, plotting the $\ell^2$ norm of solution versus coupling parameter $d$. The three branches shown correspond intersite-intersite, intersite-onsite, and onsite-onsite kink-antikinks. The inset details the intersection points of the three branches. The right panel shows six example solutions, corresponding to the labeled points on the bifurcation diagram. The insets are cartoons of the Goldstone eigenvalues for these solutions; a single marker at the origin represents a double eigenvalue at 0.}
	\label{fig:SGbifdiag}
\end{figure}

\begin{figure}
	\begin{center}
	\begin{tabular}{cc}
	\includegraphics[width=7.5cm]{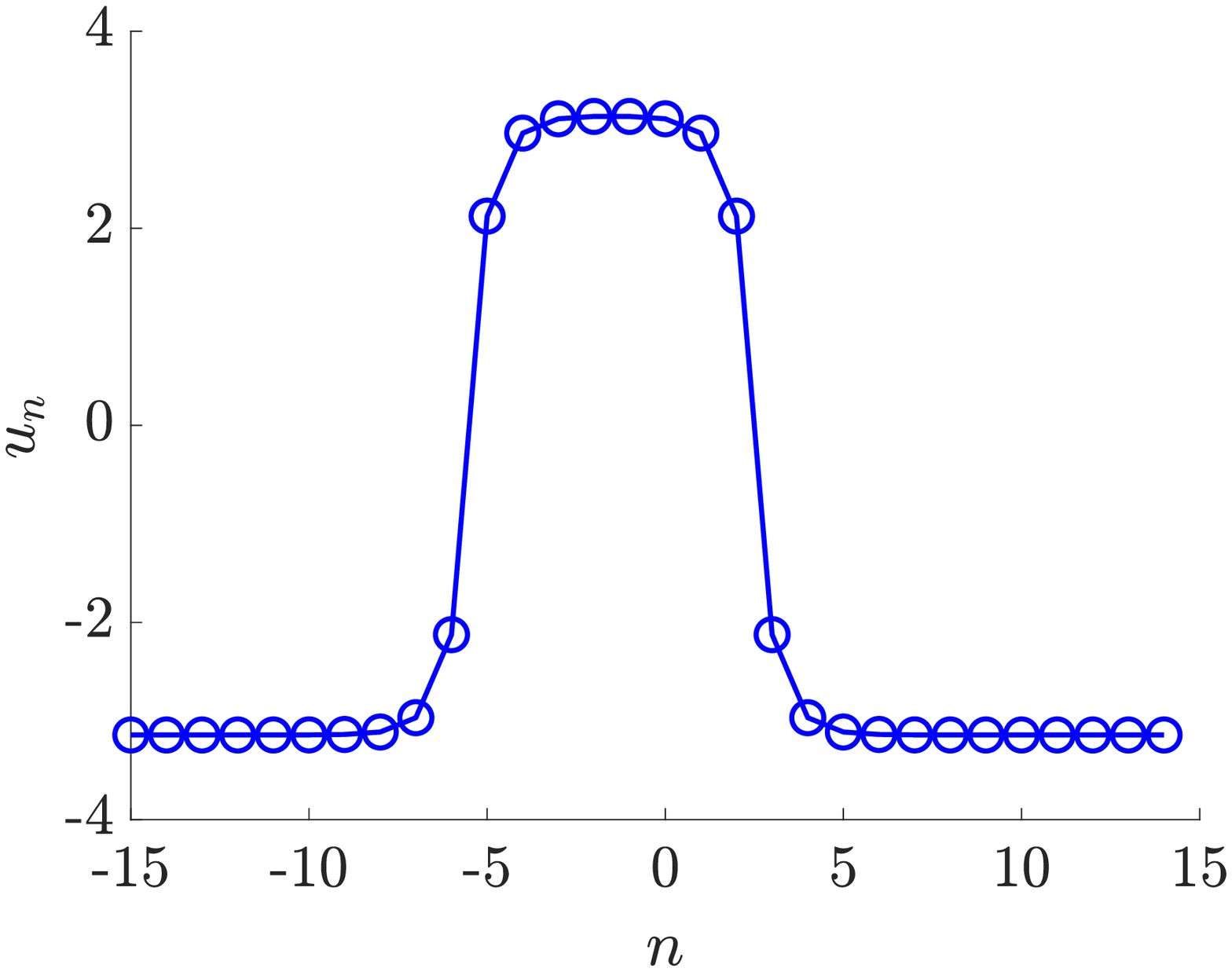} &
	\includegraphics[width=7.5cm]{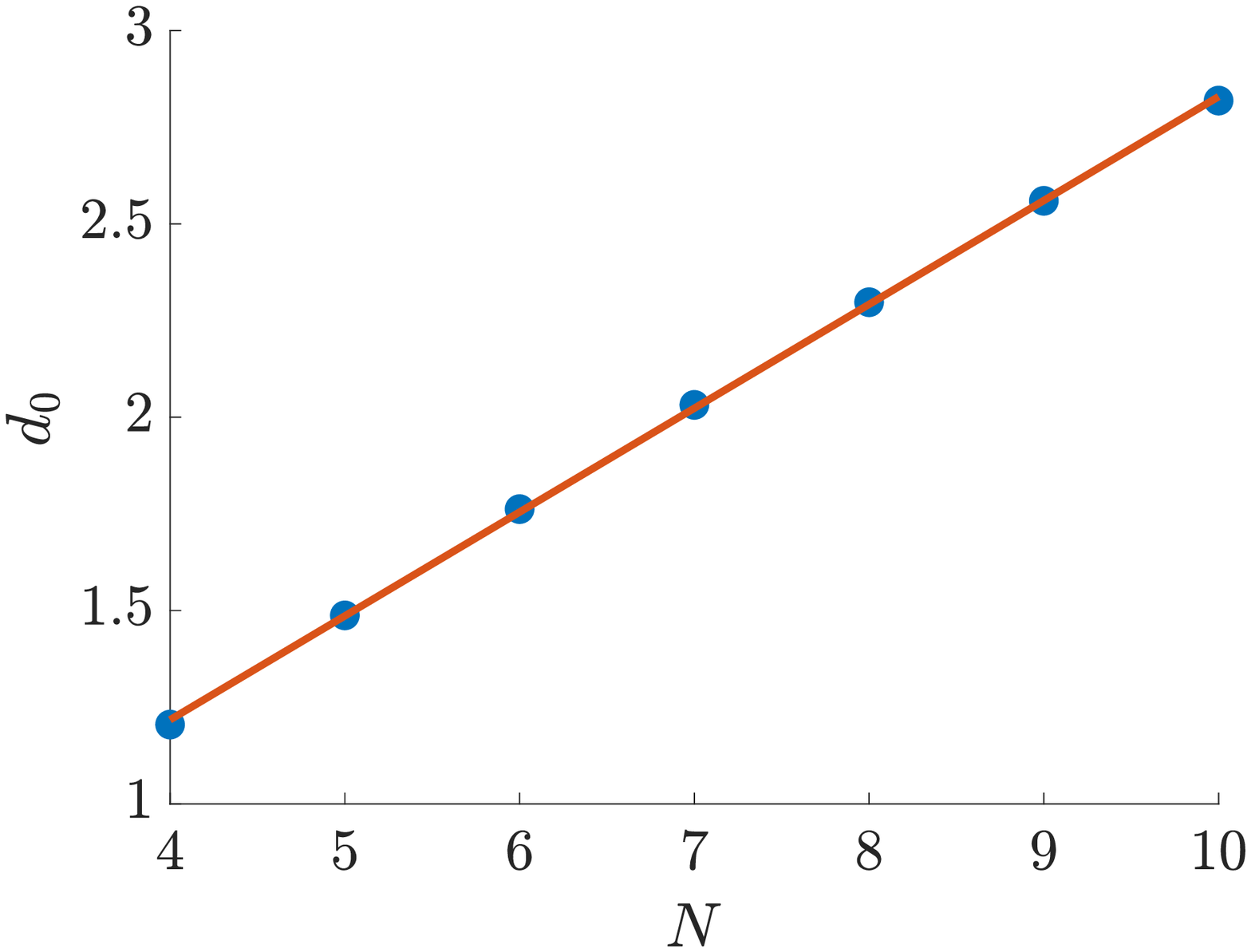}
	\end{tabular}
	\end{center}
	\caption{The left panel shows a kink-antikink solution with $N_1 = 8$ for $d = 0.5$. The right panel shows the turning point $d_0$ vs $N$ for kink-antikink waveforms together with least squares linear regression line.}
	\label{fig:kak}
\end{figure}

We construct a kink-antikink (\cref{fig:kak}, middle panel) by parameter continuation in the coupling parameter $d$ from the AC limit using the software package AUTO. As an initial condition, we use a kink-antikink composed of two intersite kinks, which at the AC limit has the form $(\cdots, -\pi, -\pi, \pi, \pi, \dots, \pi, -\pi, -\pi, \dots )$, where there are $N_1$ sites in the middle of the solution which take the value $\pi$. The bifurcation diagram for the parameter continuation is shown in \cref{fig:SGbifdiag}. Notably, there is a turning point at a critical value $d_0$ (label 4 in the inset), where the kink-antikink does not exist for $d > d_0$. 
Indeed, this is natural to expect as, in the continuum
limit of the model, the attractive interaction between the
kink and antikink~\cite{MANTON1979397} 
cannot be countered by discreteness and the associated
Peierls-Nabarro barrier, and hence such a bound, stationary
state cannot exist. 
The top branch of the bifurcation diagram in \cref{fig:kak} is a kink-antikink comprising two intersite kinks, and the bottom branch is a kink-antikink comprising two onsite kinks. 
However, a direct eigenvalue count illustrates that
these branches cannot ``collide'' with each other
at a turning point (i.e., at a saddle-center bifurcation). This
because the intersite kink state is stable, while the
onsite one contains two unstable eigenvalue pairs. Hence,
there must exist also an intermediate branch with one
unstable eigenvalue pair.
Indeed, such a middle branch exists, and is an asymmetric kink-antikink comprising one intersite kink and one onsite kink, which meets the bottom branch at a pitchfork bifurcation point 
(indeed, there are two realizations of the 
intersite-onsite kink which are mirror images of each other, as
discussed in more detail below) at a value of $d$ slightly smaller than $d_0$. The center branch of the bifurcation diagram in \cref{fig:SGbifdiag} is composed of two branches: solutions on one branch are an intersite kink followed by an onsite antikink (label 2), and solutions on the other branch are an onsite kink followed by an intersite antikink (not shown in the figure). Solutions on these two branches are left-right mirror images of each other and have the same $\ell^2$ norm at the same value of $d$. Insets in the right panel of \cref{fig:SGbifdiag} show the Goldstone eigenvalue pattern for the kink-antikink solutions. Pairs of Goldstone eigenvalues collide at the origin at the turning point and the pitchfork bifurcation points. 
First, the collision of the asymmetric intersite-onsite
branch with the onsite-onsite one takes place: as a result
of this pitchfork bifurcation, the asymmetric branch disappears, and the symmetric waveform emerging thereafter has only
a single pair of unstable eigenvalues. Then, at the turning point, it collides in turn with the intersite-intersite branch, and the branches disappear past this critical
point $d_0$ of the collision.
The parameter continuation suggests a linear relationship between the critical value of the coupling parameter $d_0$ and the separation distance $N$, which is shown in the right panel of \cref{fig:kak}.

\begin{figure}
	\begin{center}
	\begin{tabular}{ccc}
	\includegraphics[width=5cm]{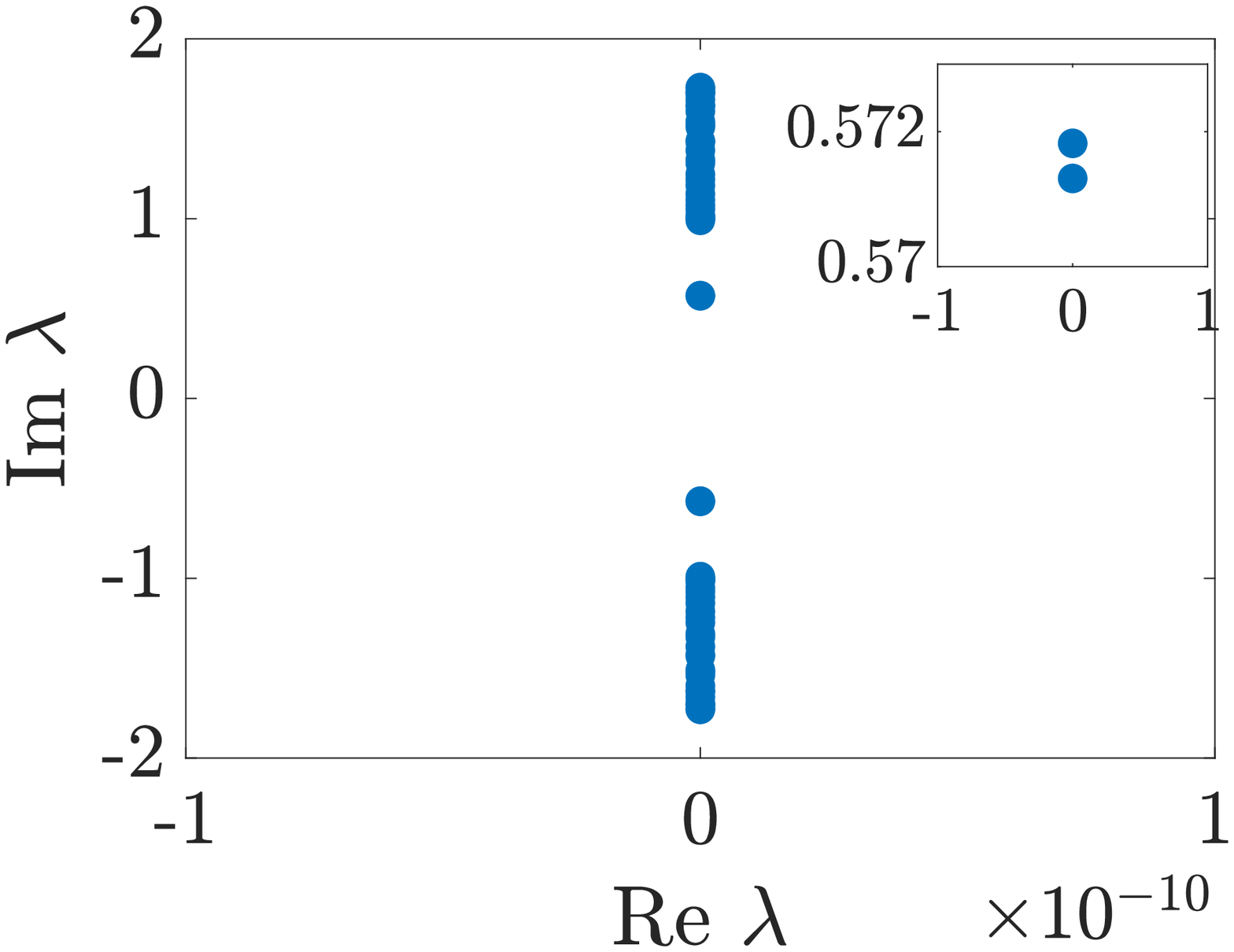}	&
	\includegraphics[width=5cm]{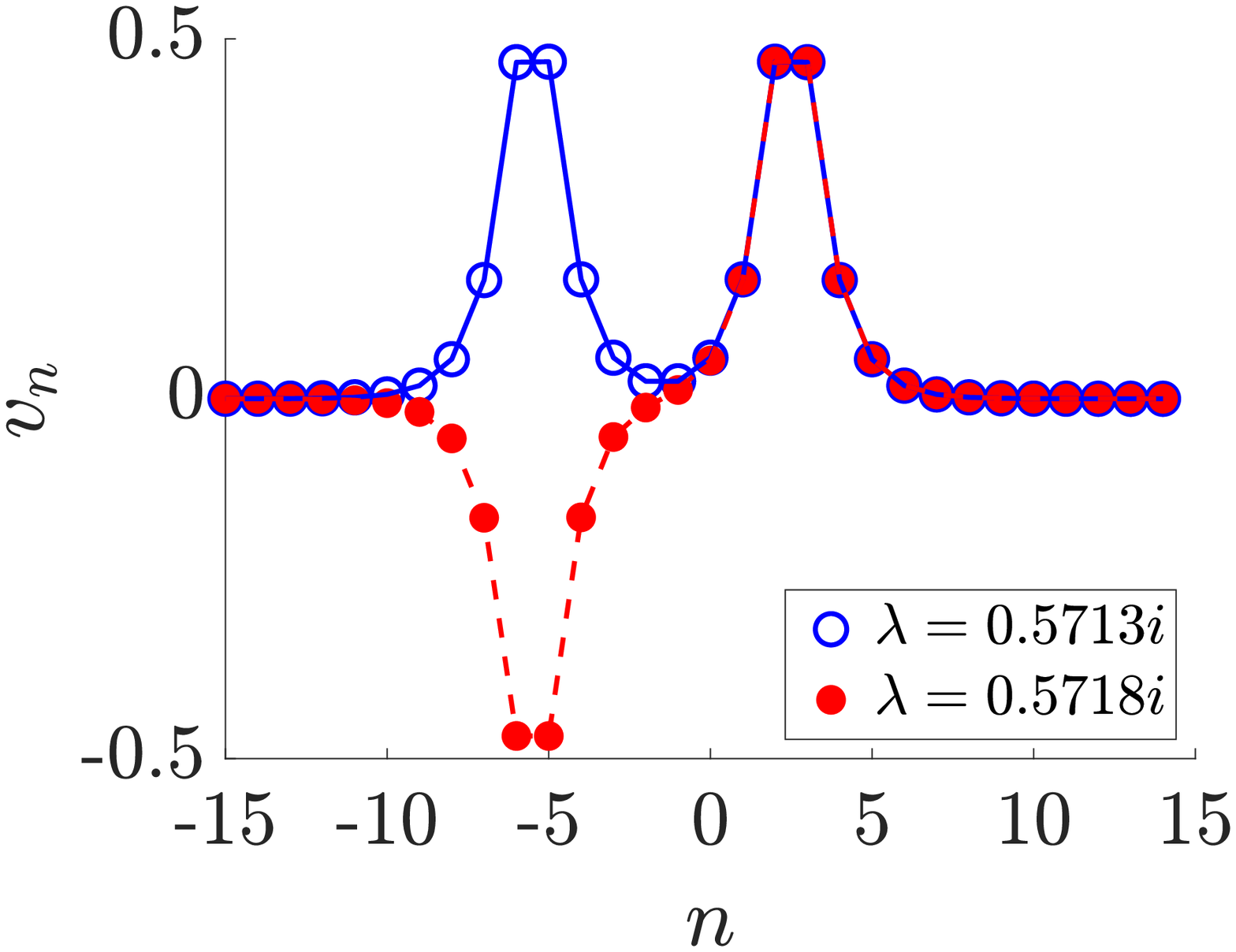} &
	\includegraphics[width=5cm]{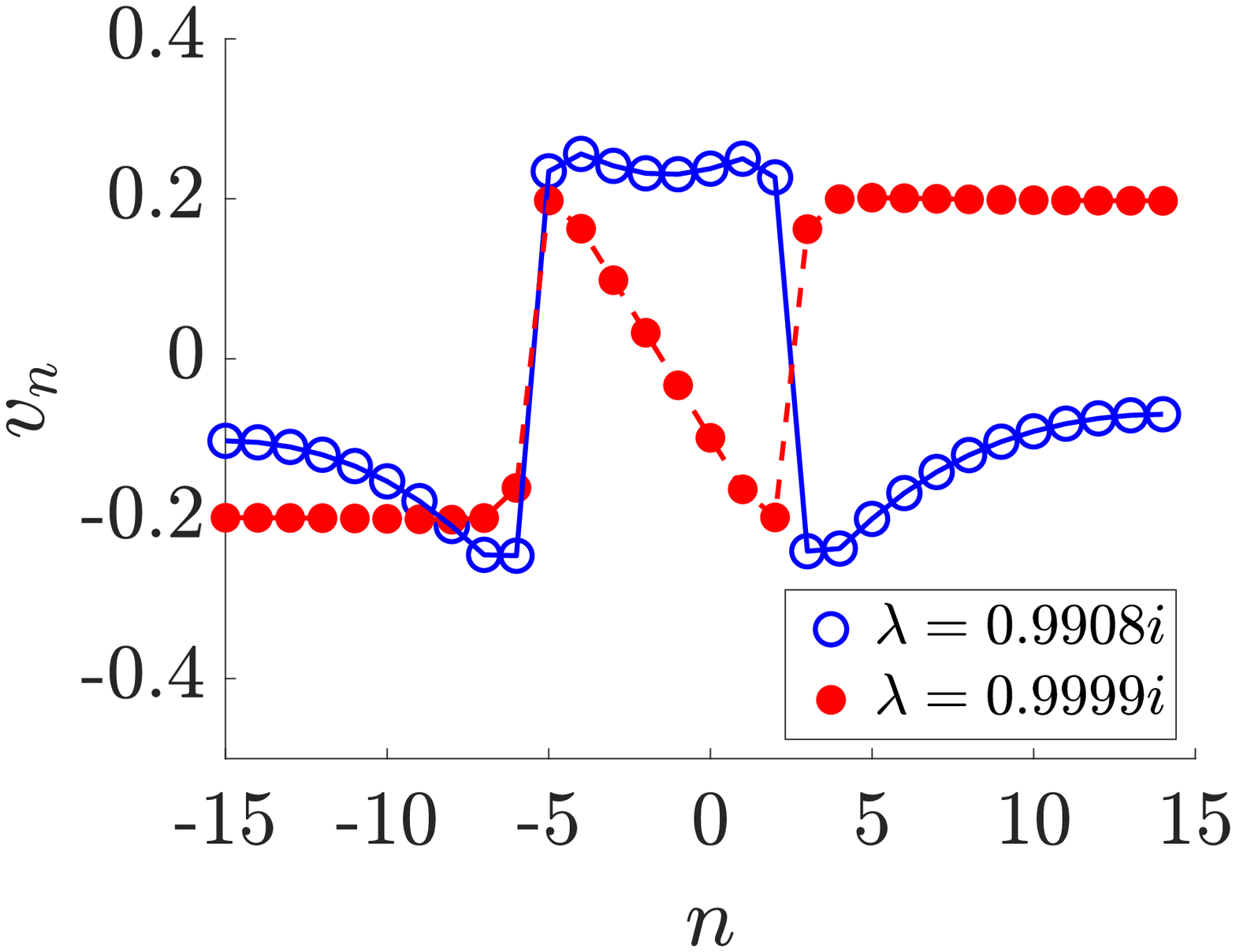}
	\end{tabular}
	\end{center}
	\caption{The left panel shows the spectrum of the kink-antikink, with an
	inset showing splitting of the Goldstone mode. The edge mode is also split (not shown). The center panel shows the eigenfunctions corresponding to the split Goldstone modes. The right panel shows the eigenfunctions corresponding to the split edge modes. $d = 0.5$, $N = 4$.}
	\label{fig:kakspec}
\end{figure}

The spectrum of the kink-antikink is similarly computed (\cref{fig:kakspec}, left panel) using MATLAB's \texttt{eig} function. As predicted by \cref{th:stability}, each element of the point spectrum splits into two eigenvalues. The eigenfunctions corresponding to the split Goldstone modes resemble two copies of the Goldstone eigenfunction of the primary kink, spliced together both in-phase and out-of-phase (\cref{fig:kakspec}, center panel). A similar phenomenon occurs with the split edge mode eigenfunctions (\cref{fig:kakspec}, right panel).

\begin{figure}
	\begin{center}
	\begin{tabular}{cc}
	\includegraphics[width=7.5cm]{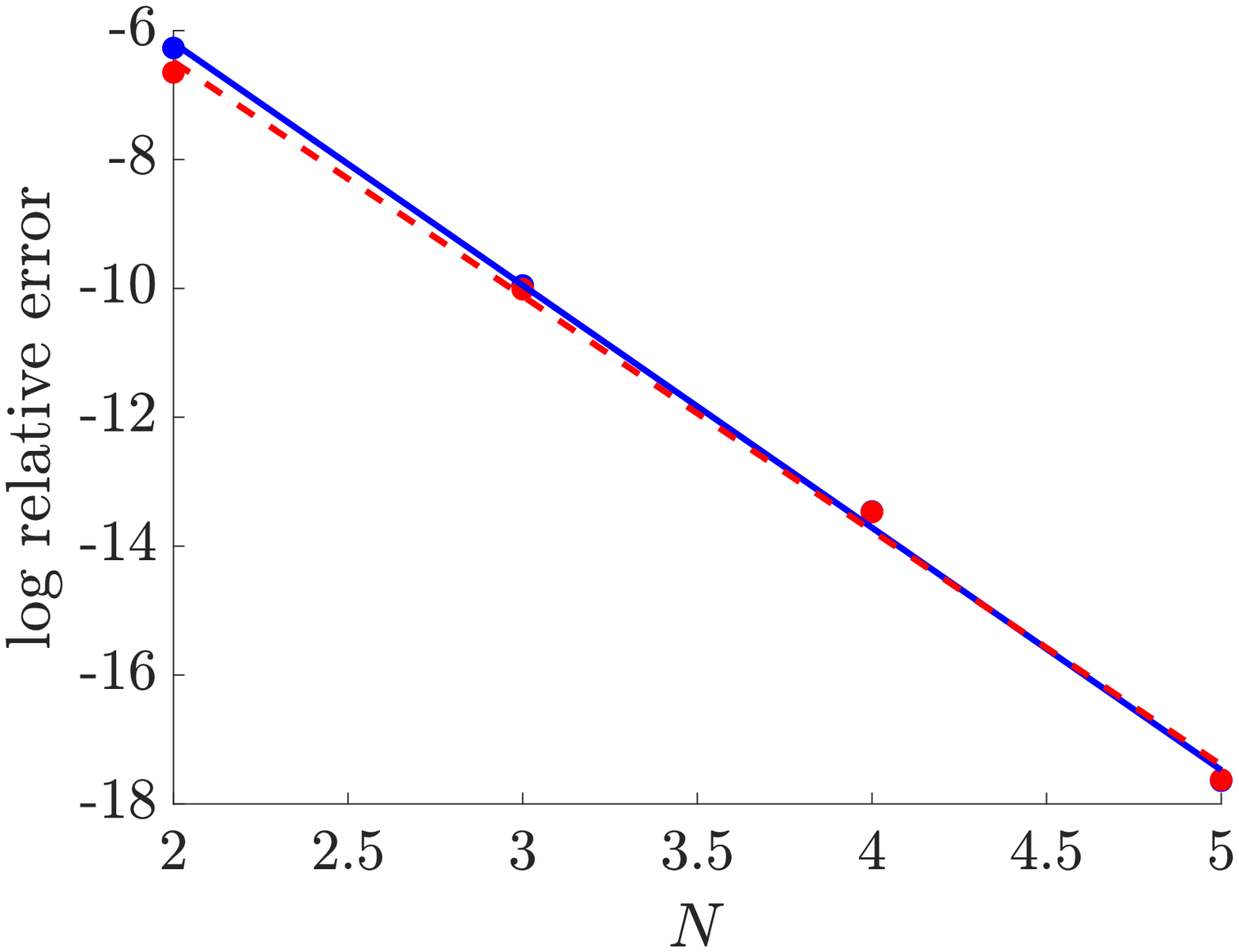} &
	\includegraphics[width=7.5cm]{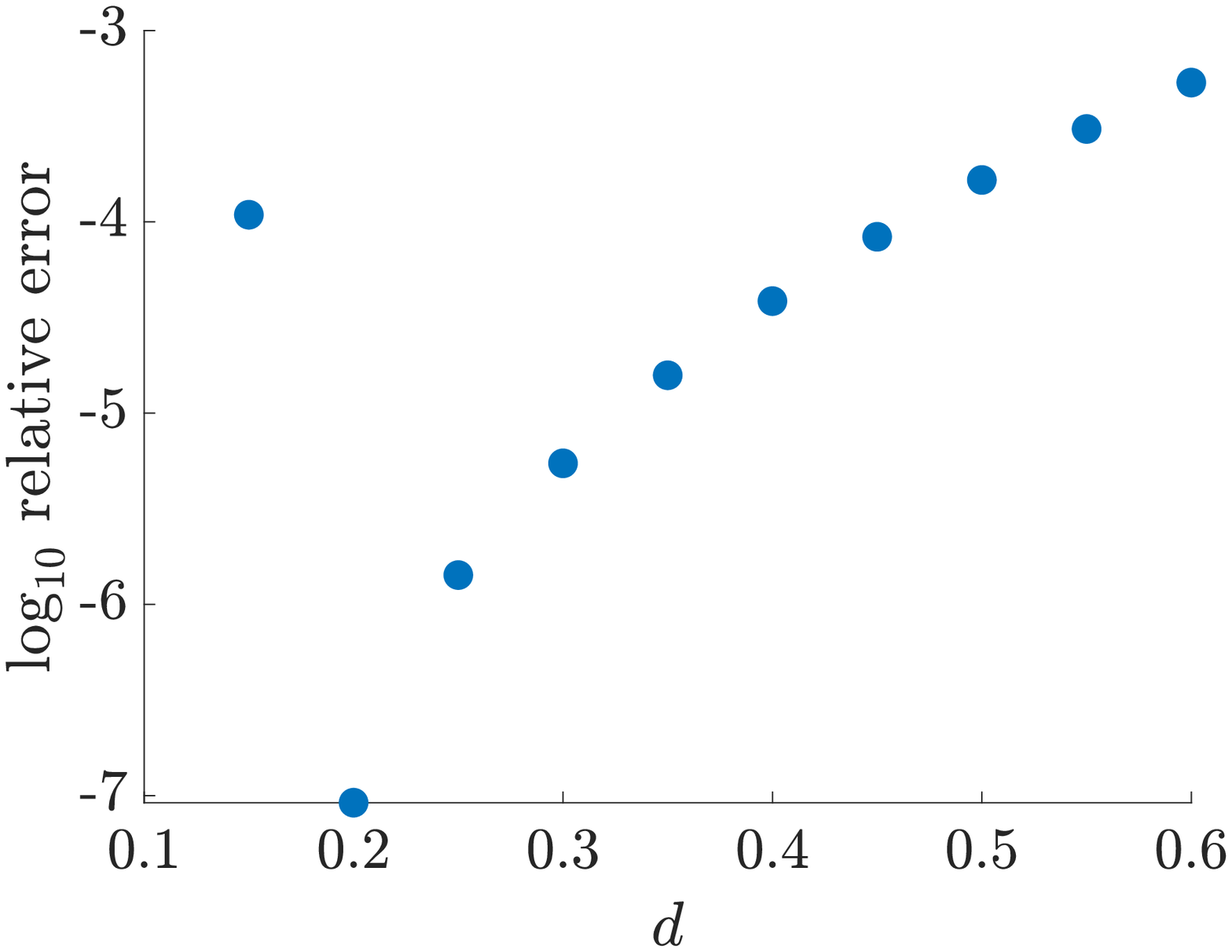}
	\end{tabular}
	\end{center}
	\caption{The left panel shows the log of the relative error in eigenvalue computation vs. $N$ for the two Goldstone eigenvalues for the kink-antikink with $d = 0.25$ together with least square linear regression lines. Blue dots and solid line denote one Goldstone eigenvalue, while red dots and dashed line the other Goldstone eigenvalue. The right panel shows $log_{10}$ of the relative error in the eigenvalue computation vs. $d$ for the two Goldstone eigenvalues for a kink-antikink waveform
	with $N = 4$.} 
	\label{fig:kakeigerror}
\end{figure}

We can compare the eigenvalues obtained from numerical computation with those predicted by \cref{corr:m23}. The left panel of \cref{fig:kakeigerror} plots the log of the relative error $| \lambda_{\text{true}} - \lambda_{\text{predicted}} |/| \lambda_{\text{true}} |$
of the two Goldstone eigenvalues vs. the separation distance $N$; the value for $\lambda$ computed with MATLAB from the linearization
around the numerically computed kink-antikink solution is used as the true value of $\lambda$.
The slope of the least square linear regression line suggests that this error is order $\mathcal{O}(r_0^{-2N})$. 
The right panel of \cref{fig:kakeigerror} plots the log of the relative error of the two Goldstone eigenvalues vs. the coupling parameter $d$. For intermediate values of $d$, the relative error is less than $10^{-3}$. The error is a minimum for approximately $d = 0.2$, and increases with increasing $d$ and as the
continuum limit is approached. (See \cite{Parker2020}*{Figure 4} for a similar phenomenon which occurs in the error plot for eigenvalues associated with double pulses in DNLS). Since the results of the \cref{th:stability} are not uniform in $d$, i.e. they hold for sufficiently large $N$ once $d$ has been chosen, we 
indeed expect a growth of the relative error for large $d$.

\begin{figure}
	\begin{center}
	\begin{tabular}{cc}
	\includegraphics[width=7.5cm]{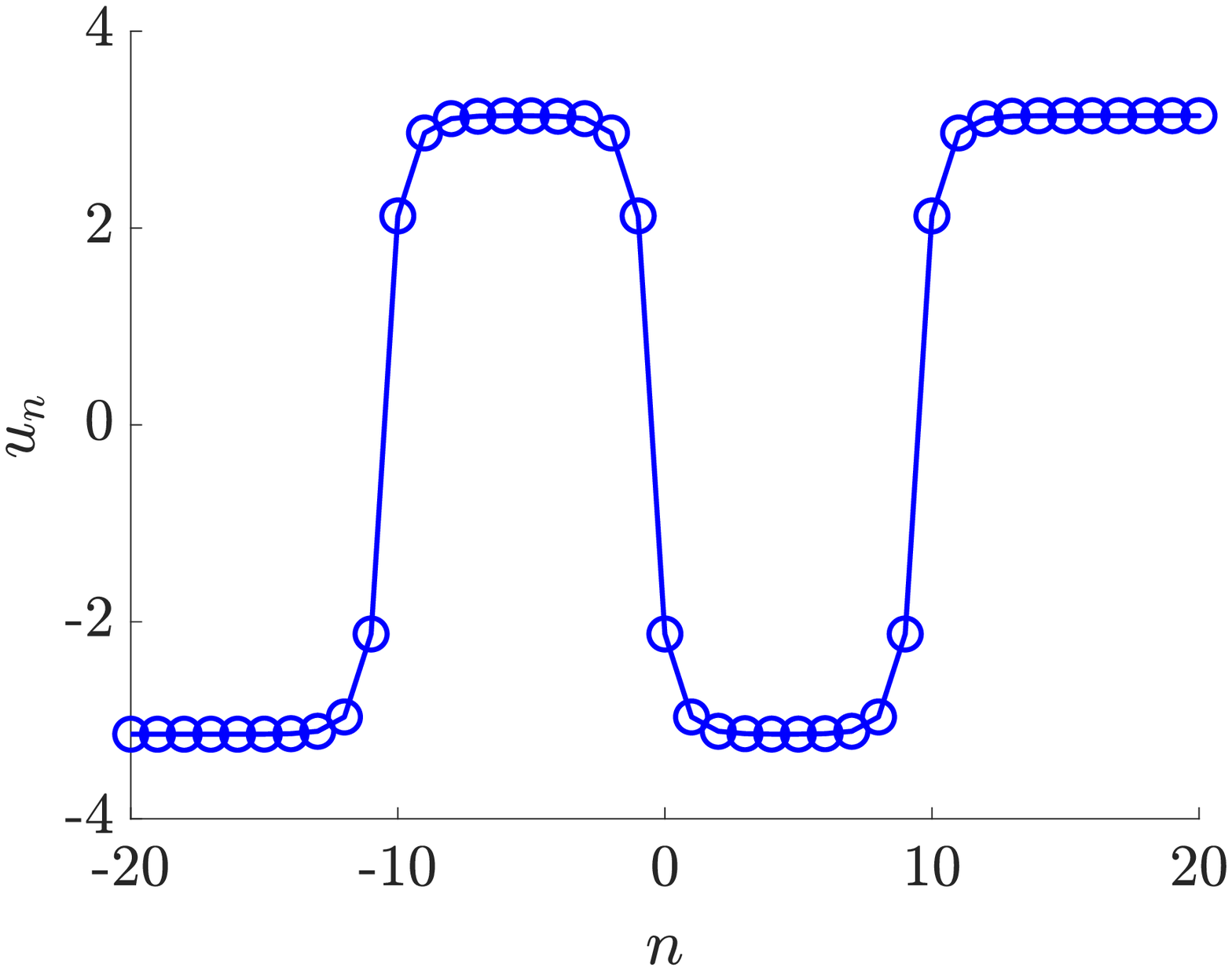} &
	\includegraphics[width=7.5cm]{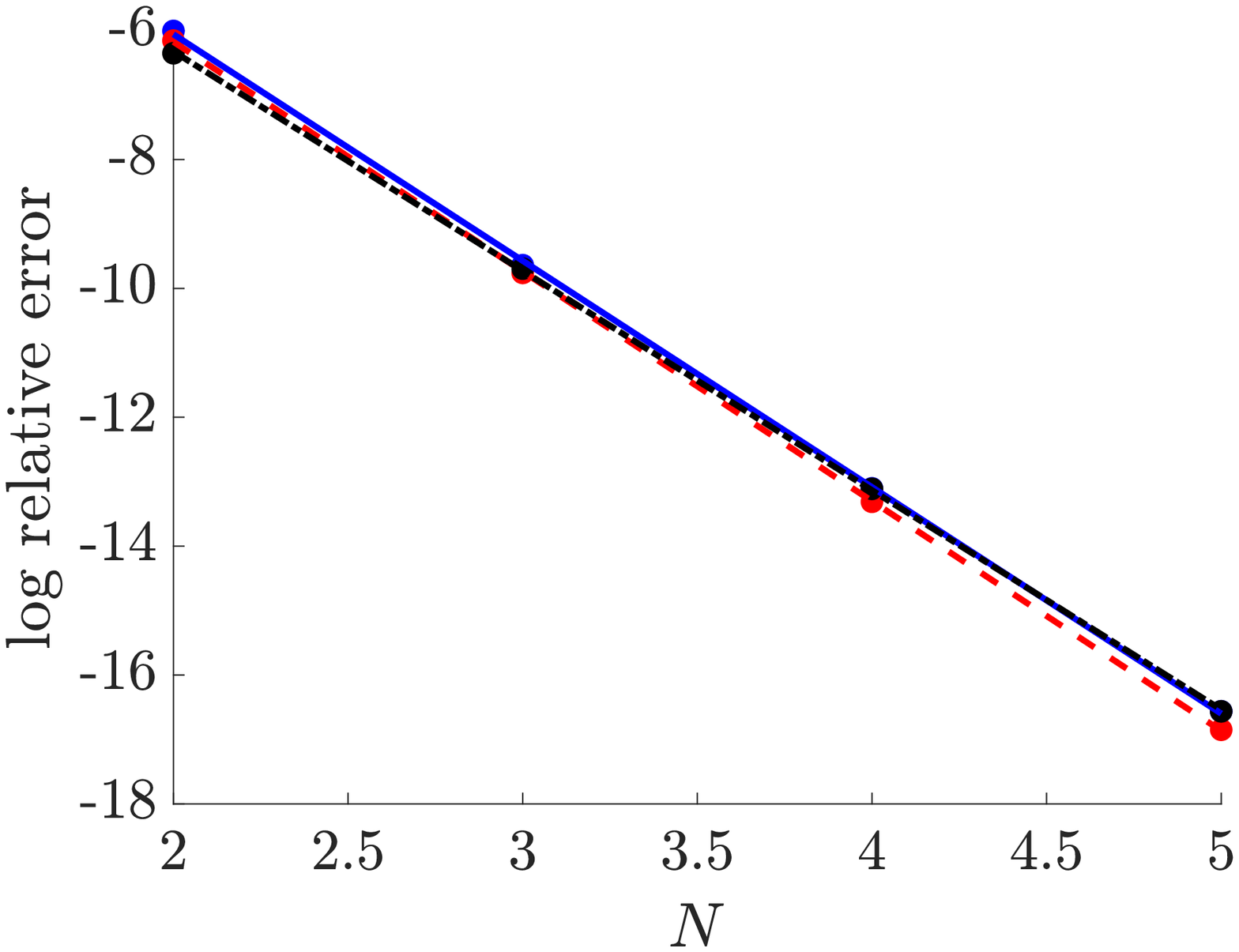}
	\end{tabular}
	\end{center}
	\caption{The left panel shows a 3-component multi-kink (kink-antikink-kink) with $N_1 = N_2 = 8$ and $d = 0.25$. 
	The right panel shows the log of the relative error in eigenvalue computation vs. $N$ for the three Goldstone eigenvalues of 3-wave multi-kink with $N_1 = N_2 = 2N$ and $d = 0.25$. The lines are least square linear regression fits. Each line among the 3 very proximal solid, dashed and dash-dotted
	lines corresponds to one of the three Goldstone eigenvalues.} 
	\label{fig:3p}
\end{figure}

We can obtain similar results for higher order multi-kinks. An example of a three-component multi-kink is shown in the left panel of \cref{fig:3p}. We can again compare the eigenvalues obtained from numerical computation with those predicted by \cref{corr:m23}. The right panel of \cref{fig:3p} shows the relative error of the eigenvalue computations for the three Goldstone eigenvalues.
One can again observe the particularly good agreement
of the theory and the computation, especially so for
larger values of $N$.

\begin{figure}
	\begin{center}
	\begin{tabular}{cc}
	\includegraphics[width=7.5cm]{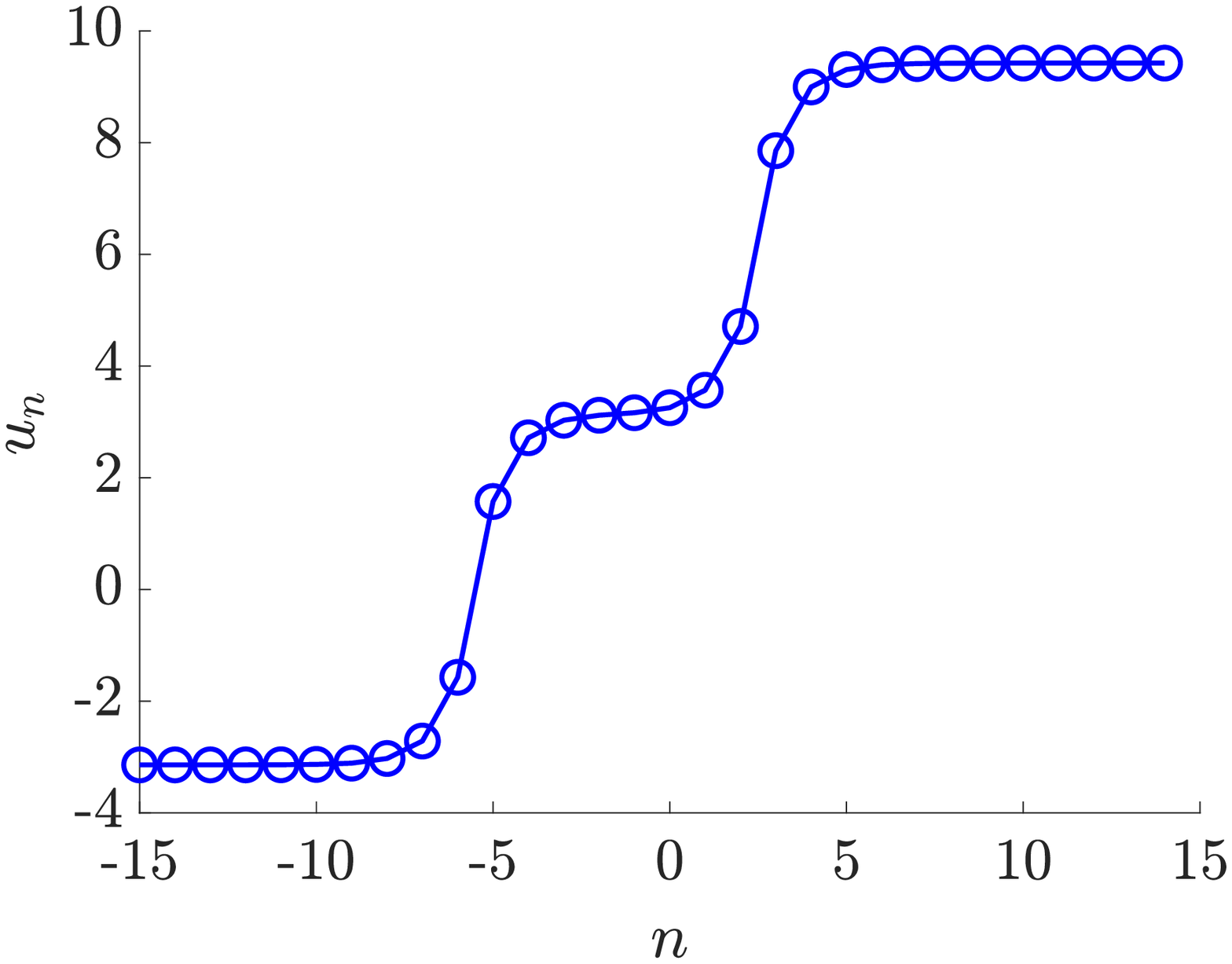} &
	\includegraphics[width=7.5cm]{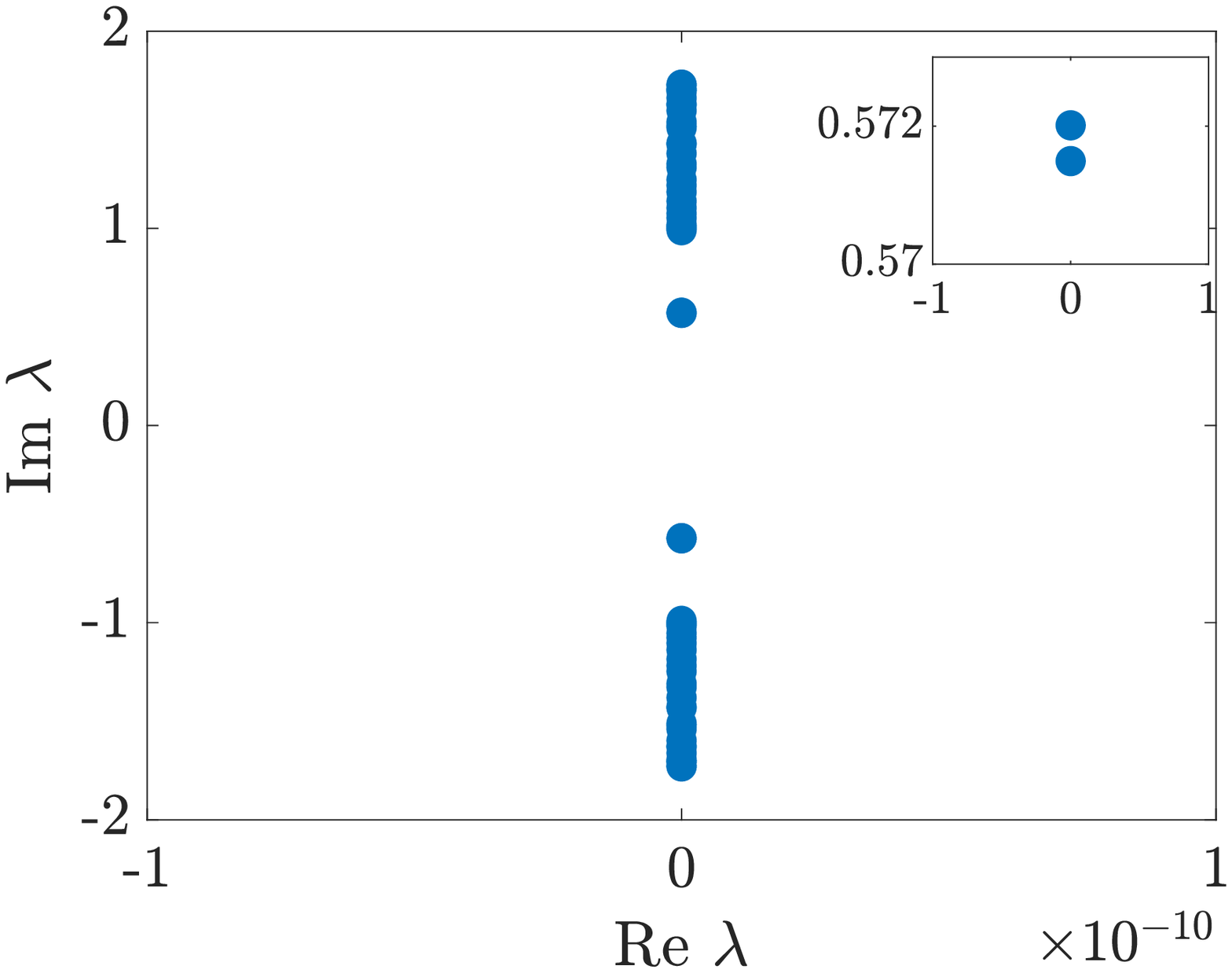}
	\end{tabular}
	\end{center}
	\caption{The left panel shows a kink-kink solution with $N_1 = 8$. The right panel shows the spectrum of the kink-kink waveform, with the inset showing splitting of Goldstone mode. The edge mode is also similarly 
	split (not shown). $d = 0.25$. }
	\label{fig:kk50}
\end{figure}

In addition, for the sine-Gordon equation, we can have generalized multi-kink solutions, in which each kink or antikink in the sequence connects two adjacent saddle equilbria (see \cref{rem:SGmulitkinks}). An example of a kink-kink is shown in \cref{fig:kk50}. The spectrum is almost identical to that of the kink-antikink with the same parameters in \cref{fig:kakspec}.

Finally, we perform timestepping simulations of the multikink
states. The topic of the evolution of
a single (primary) kink was touched upon in the work
of~\cite{KevrekidisWeinstein2000}, where dynamical
evolution experiments were performed in the case
of a single kink with the Goldstone and/or edge
modes excited. It was found that a mechanism of
resonance of the point spectrum mode harmonics
(i.e., $2 \lambda$, $3 \lambda$ etc., 
where $\lambda$ is the Goldstone
or edge eigenvalue) with the continuous spectrum
led to nonlinearity-induced power law decay of
relevant mode amplitudes. Here, we instead focus on the evolution
of the central nodes of the kinks (e.g. nodes
$-1$ and $0$ for a kink located at the origin) for the
realm of multi-wave, kink-antikink structures.
For a timestepping scheme, since the spatial component is already discretized, we use a symplectic and symmetric implicit Runge-Kutta method \cite{HairerBook}, as suggested in \cite[Section 2.5]{Duncan1997}, to preserve the symplectic structure of the Hamiltonian equation \cref{eq:dSG}. Specifically, we use the MATLAB implementation of the \texttt{irk2} scheme of order 12 from \cite{Hairer2003}. For boundary conditions, we use the discrete analogue of Neumann boundary conditions.

\begin{figure}
	\begin{center}
	\begin{tabular}{cc}
	\includegraphics[width=7cm]{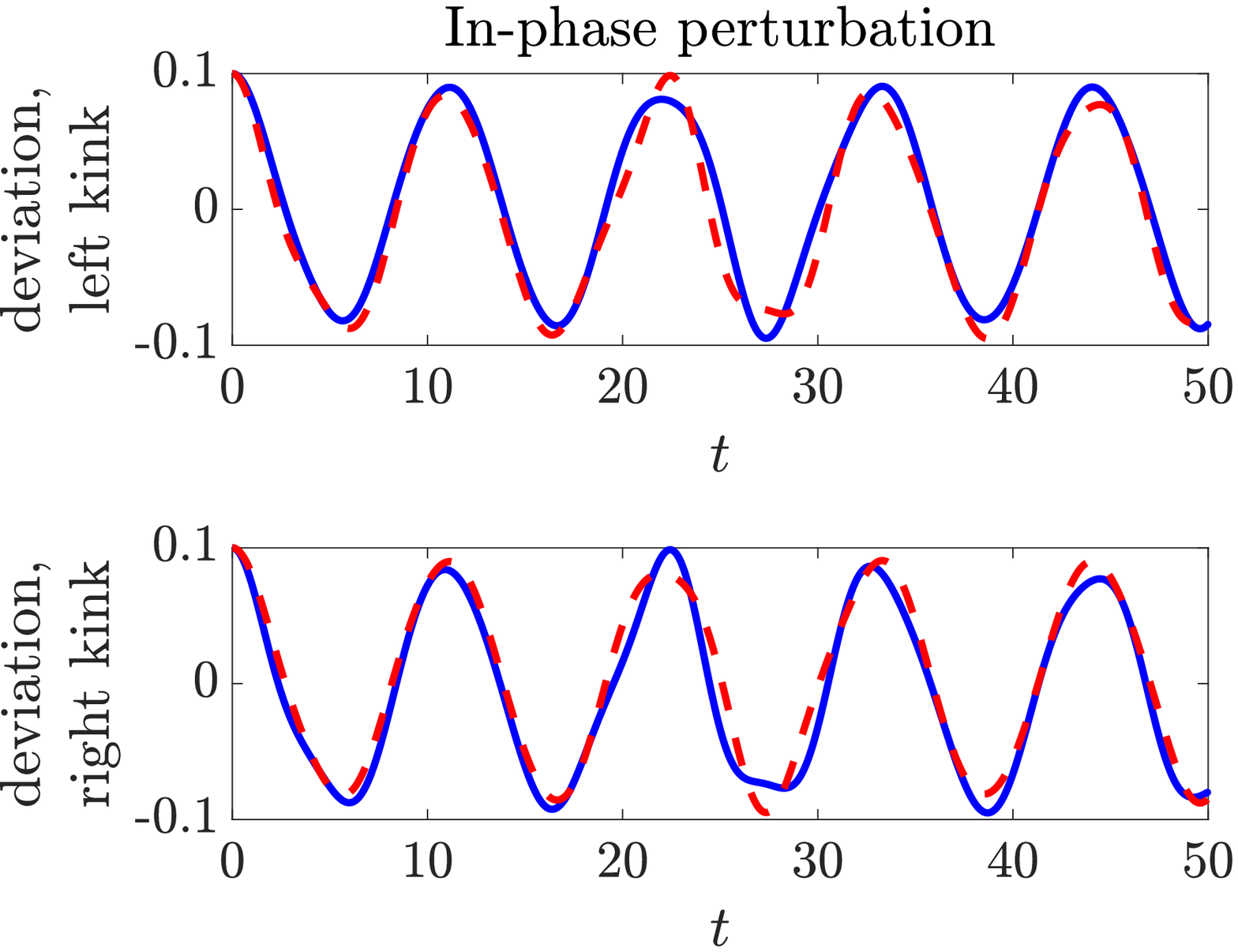} &
	\includegraphics[width=7cm]{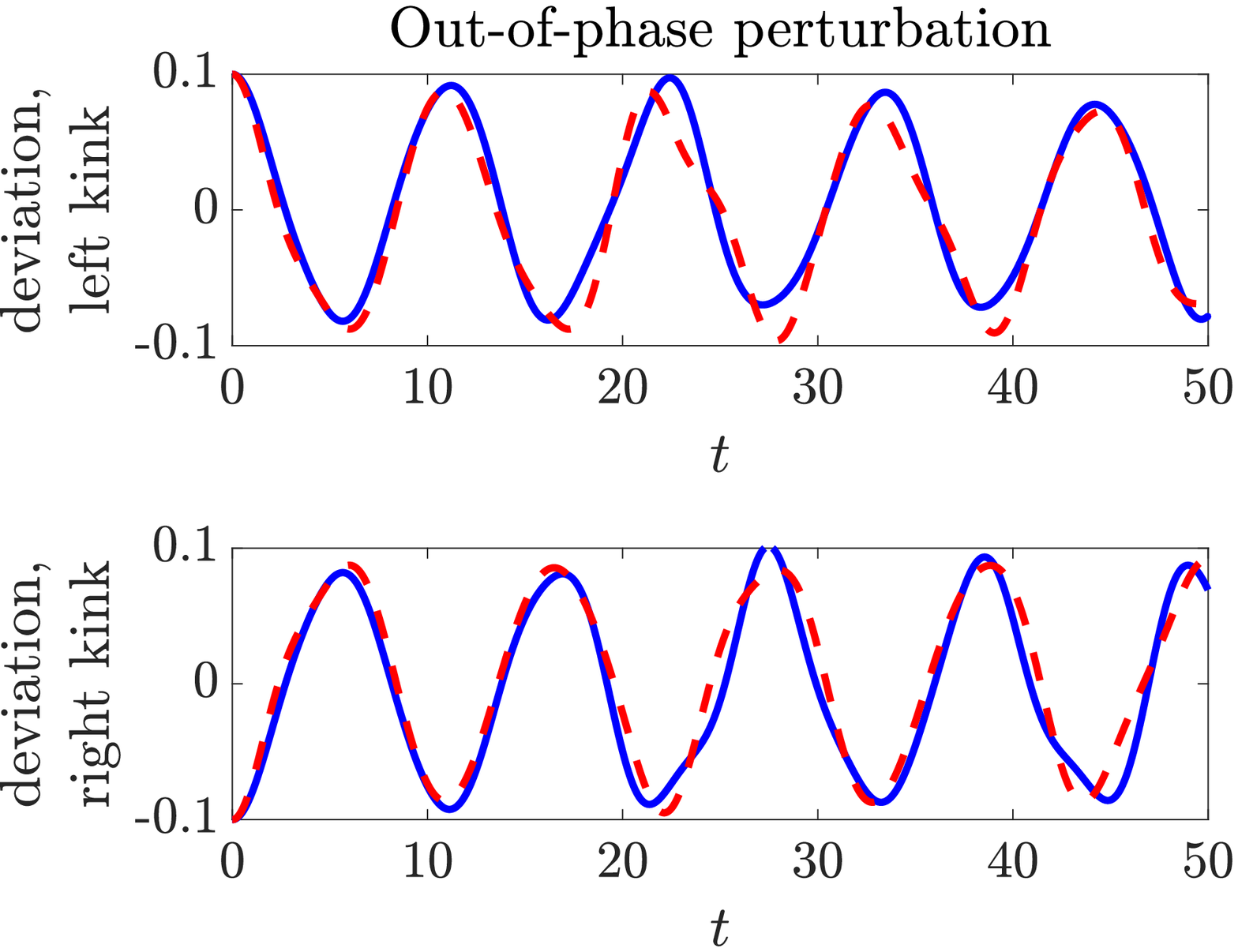} \\
	\includegraphics[width=7cm]{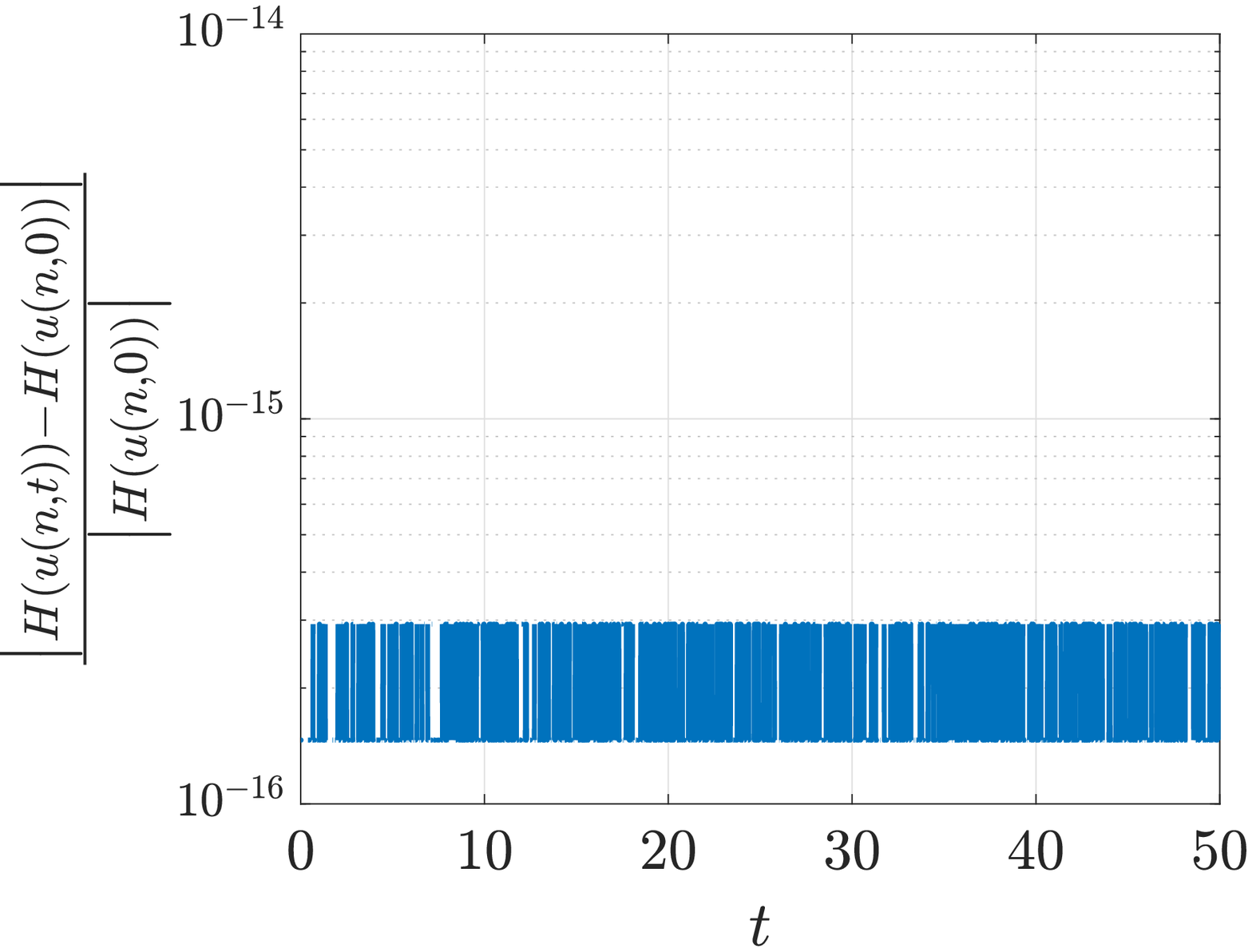} &
	\includegraphics[width=7cm]{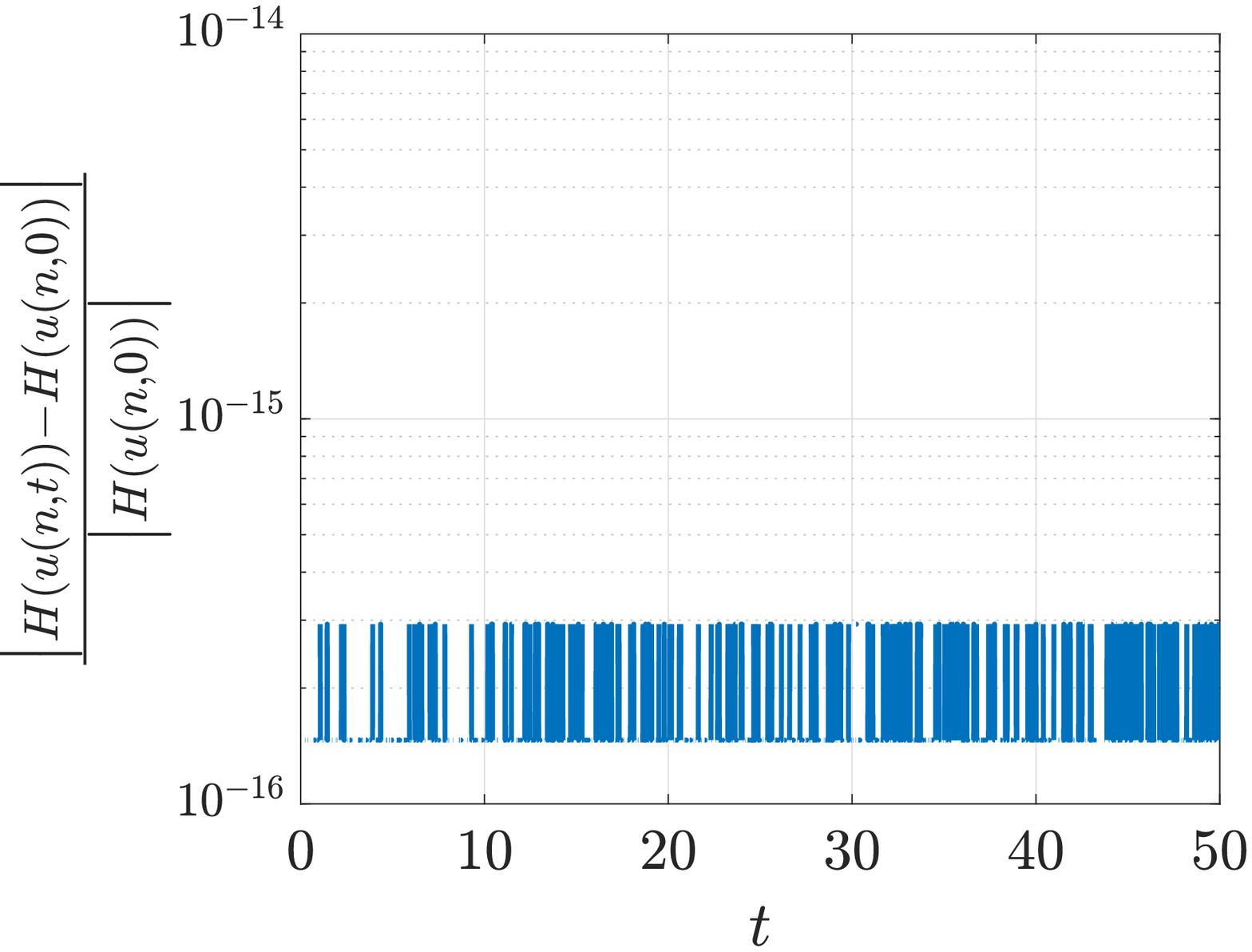} 
	\end{tabular}
	\end{center}
	\caption{Timestepping of perturbations for a kink-antikink solution. The left column shows an in-phase perturbation, where the central nodes of both kinks are perturbed by 0.1 in the same direction. The right column shows an out-of-phase perturbation, where the central nodes of the left kink are perturbed by 0.1, and the central nodes of the right kink are perturbed by -0.1. First and second rows plot deviations from the stationary kink solution of the two central nodes for the left and right kinks (respectively). The third row is a semilog plot of the relative difference of the energy $H(u(n, t))$ from the starting energy $H(u(n, 0))$. $d = 0.5$, $N = 400$ grid points, timestepping using symplectic and symmetric implicit Runge-Kutta method \texttt{irk2} with step size 0.01.} 
	\label{fig:kaktimestep}
\end{figure}

For a kink-antikink constructed from two intersite kinks, the pair of Goldstone mode eigenfunctions suggests that there will be two corresponding normal modes of oscillations for the central nodes of the two kinks: an in-phase mode and an out-of-phase mode. These can be seen in the left and right columns of \cref{fig:kaktimestep}. In addition, we plot the energy of the solution $H(u(n, t))$ as $t$ evolves (\cref{fig:kaktimestep}, bottom row). The relative deviation of the energy from its initial value is less than $10^{-15}$ over the time interval of the simulation.
There will similarly be two normal modes corresponding to the pair of edge mode eigenfunctions.
These panels are representative of the possible in- and out-of-phase motion of the multiple 
coherent structures.

\begin{figure}
	\begin{center}
	\begin{tabular}{cc}
	\includegraphics[width=7cm]{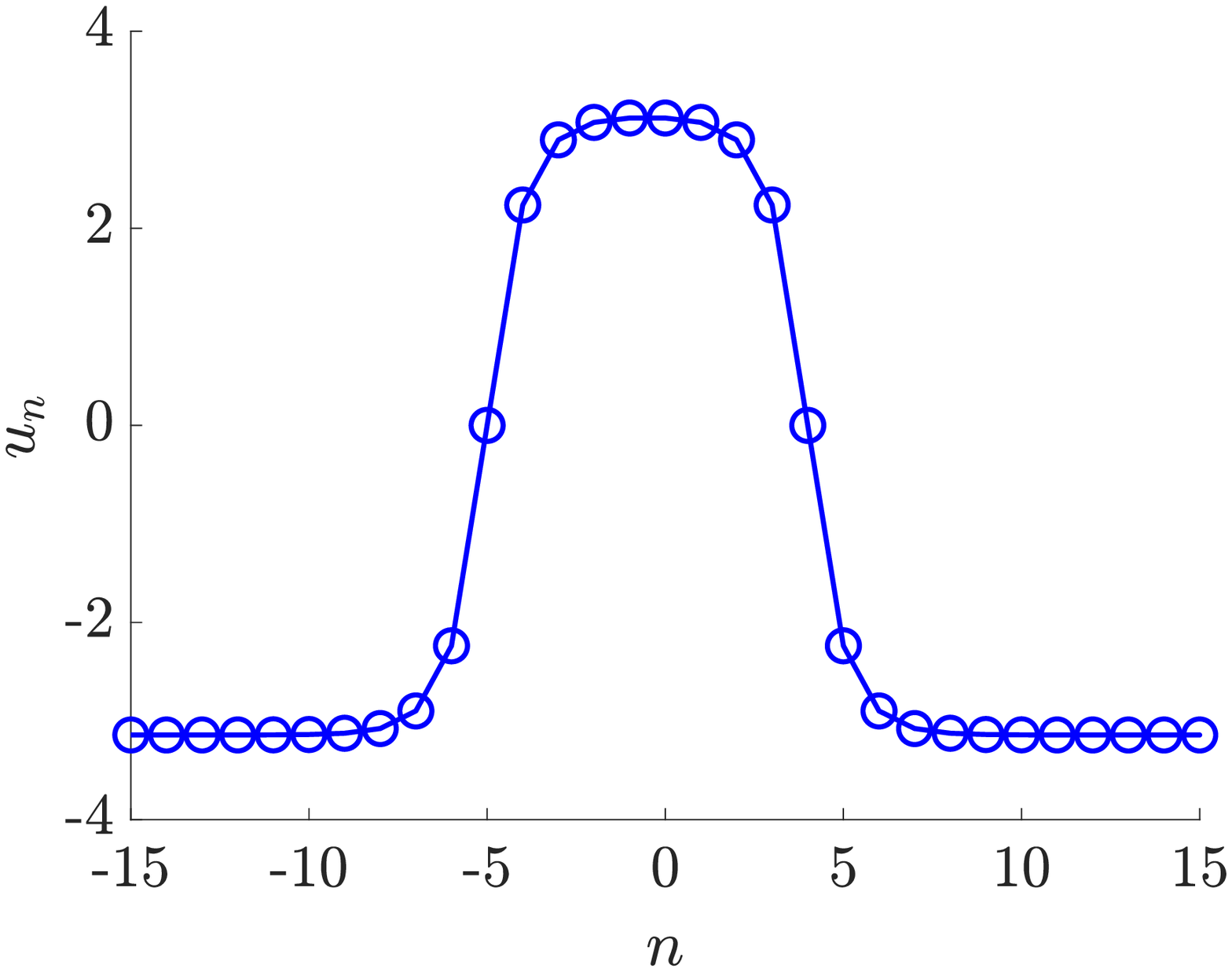} &
	\includegraphics[width=7cm]{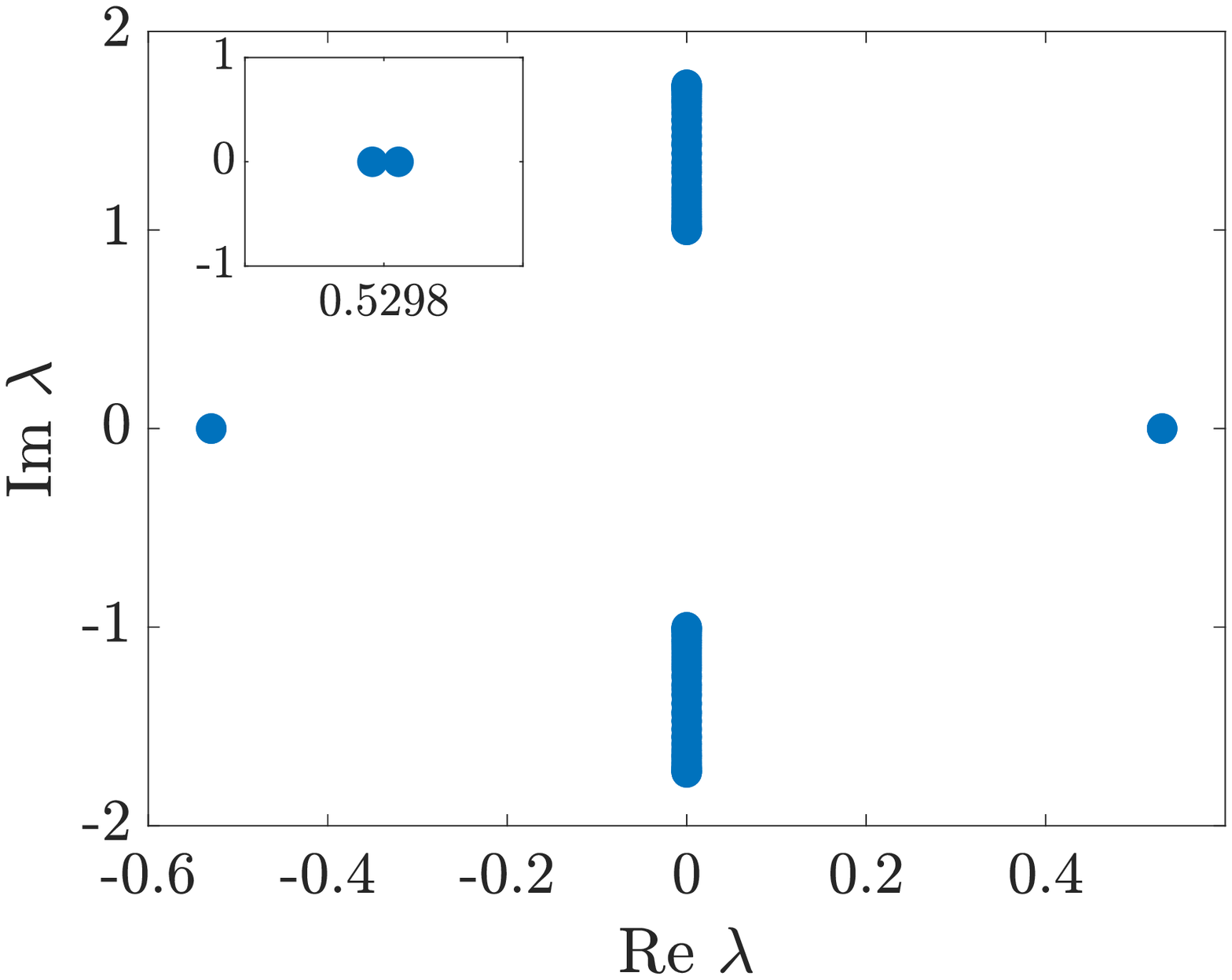} \\
	\includegraphics[width=7cm]{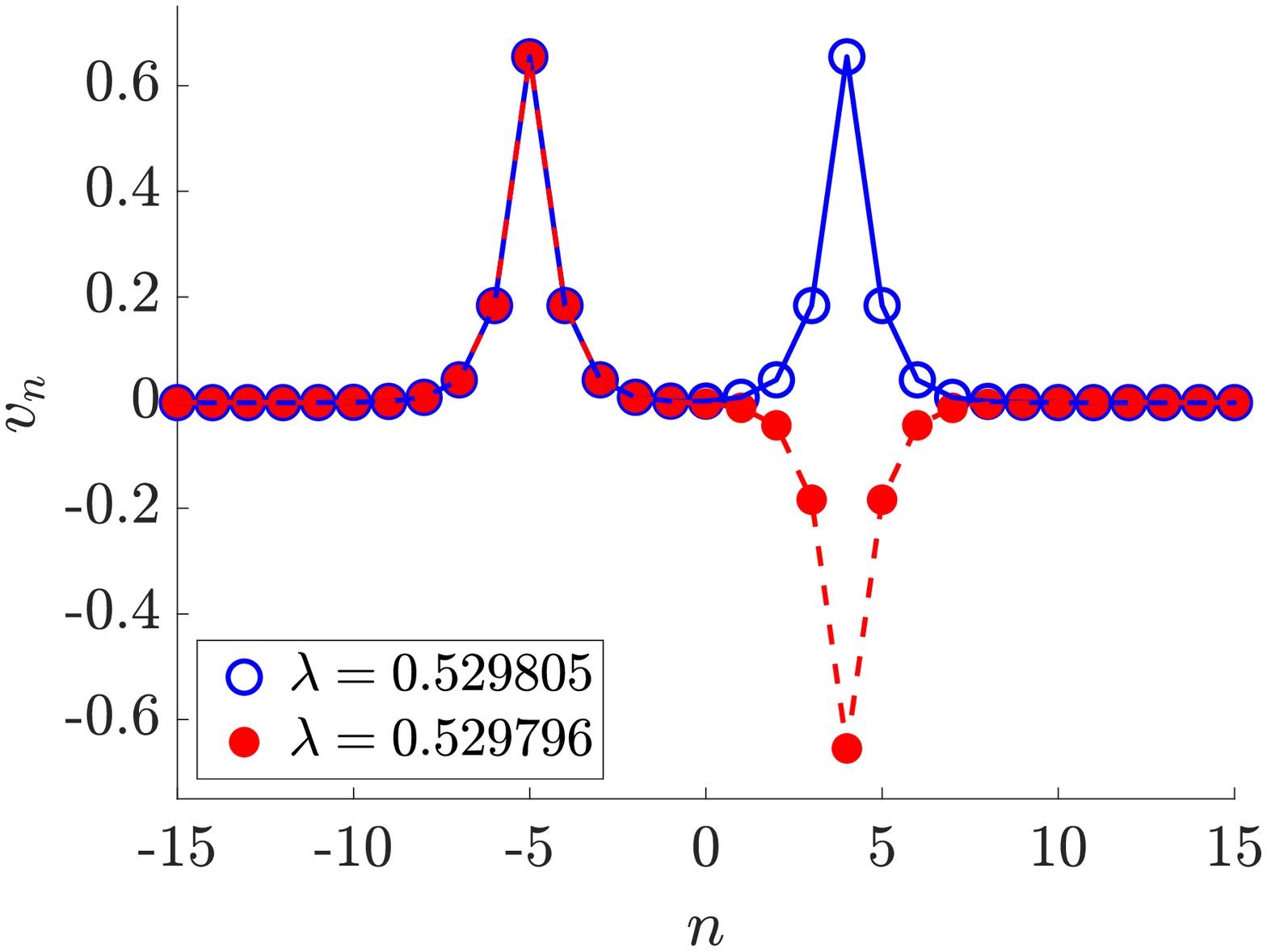} &
	\includegraphics[width=7cm]{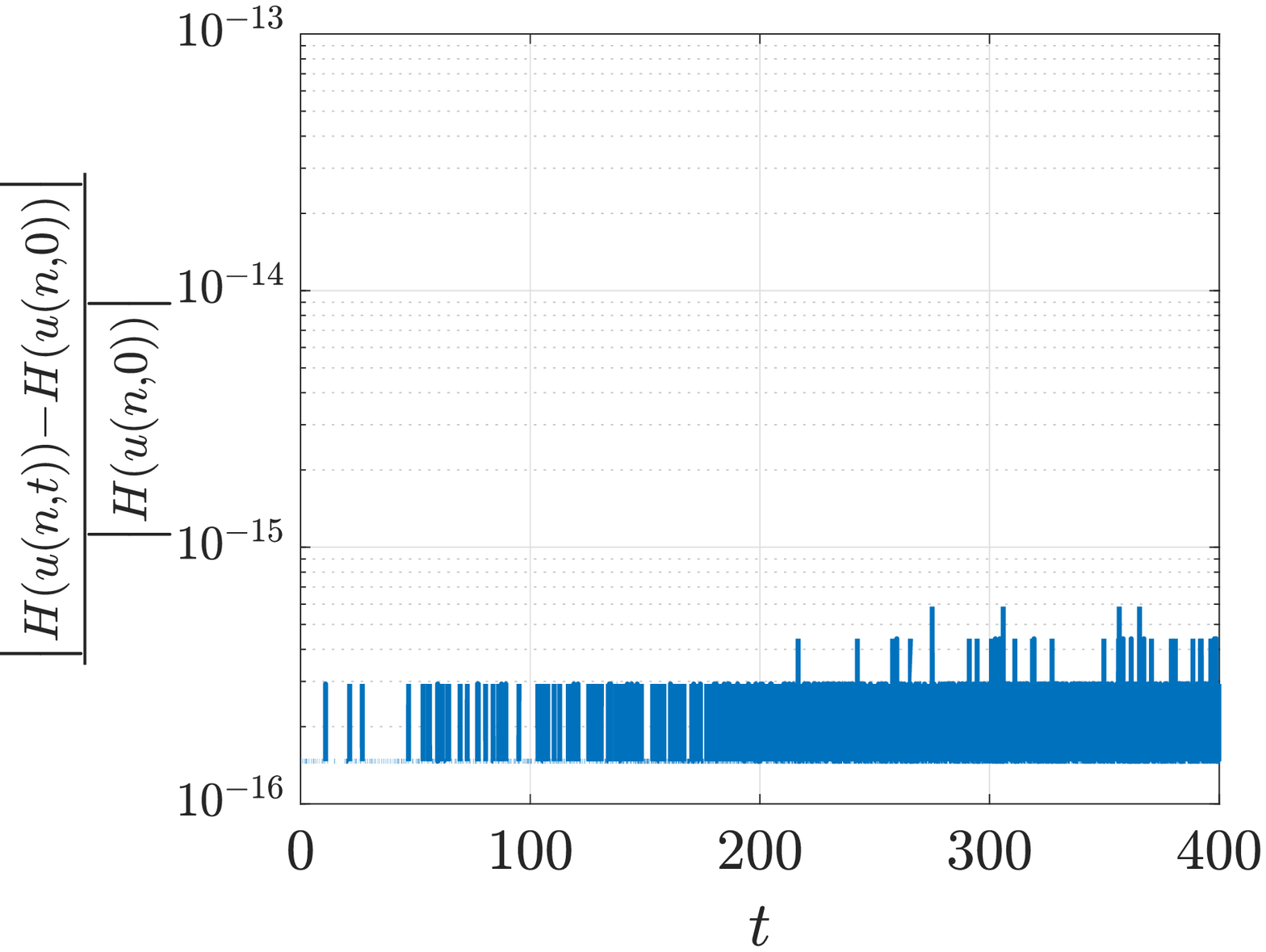} \\
	\includegraphics[width=7cm]{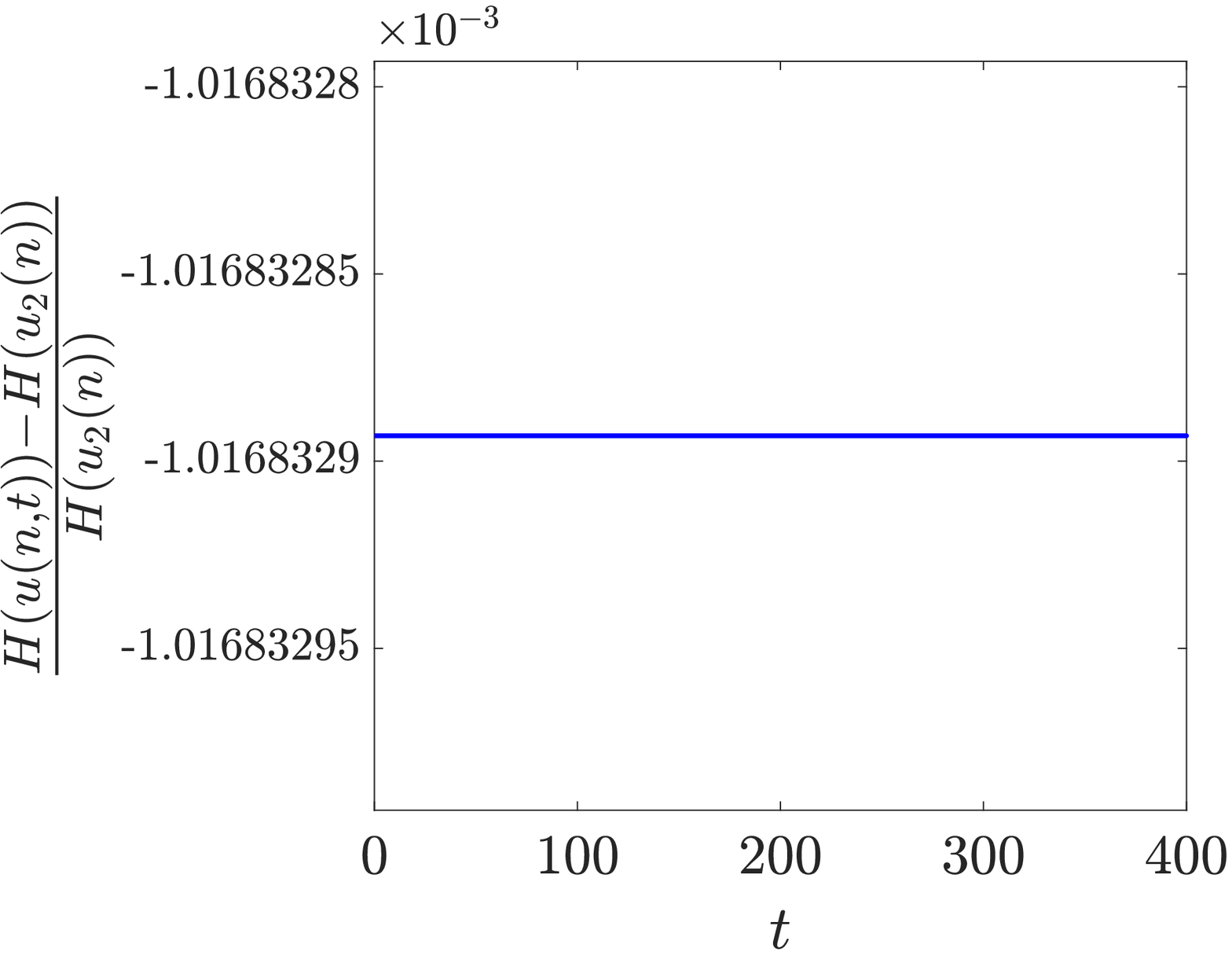} &
	\includegraphics[width=7cm]{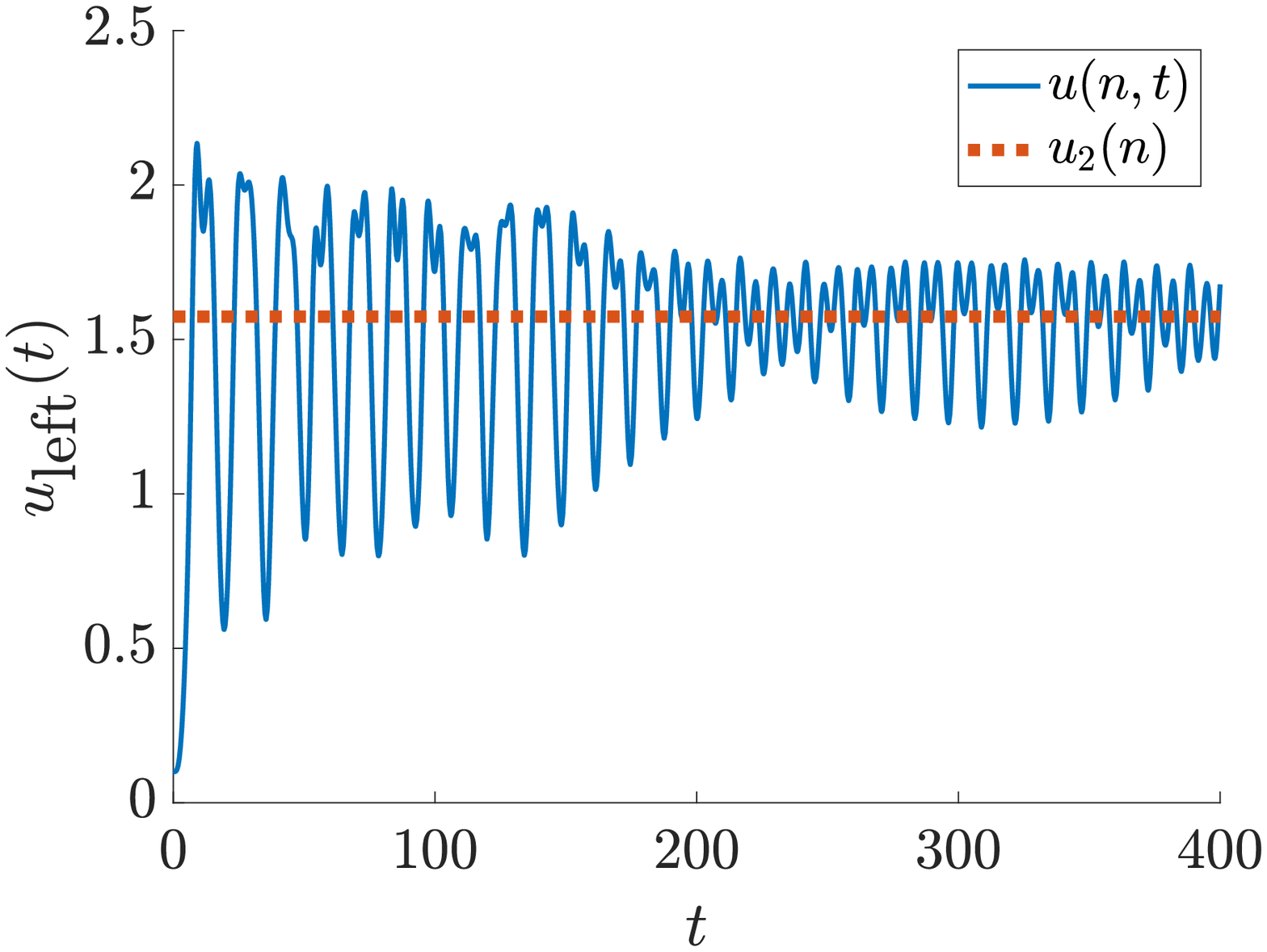}
	\end{tabular}
	\end{center}
	\caption{Kink-antikink solution constructed from two onsite kinks with $N_1 = 8$ (top left). Spectrum of onsite kink-antikink bearing two unstable near-identical Goldstone modes (top right); see also the inset discerning between the two modes. Eigenfunctions corresponding to the split Goldstone modes (middle left). Remaining plots show time evolution of perturbation $u(n,t)$ of onsite kink-antikink; initial condition obtained by adding 0.1 to the two central nodes with $u_n = 0$, leading to the destabilization of the structure. Middle right is semilog plot of the relative difference of the energy $H(u(n, t))$ from the starting energy $H(u(n, 0))$. Bottom left is energy difference between the perturbed onsite kink-antikink $u(n,t)$ and the neutrally stable static intersite kink-antikink $u_2(n)$. Bottom right shows left central node of perturbed onsite kink-antikink $u(n,t)$ (blue solid line) oscillating about left central node of static inter1site kink-antikink $u_2(n)$ (horizontal, dotted orange line). $N=400$ grid points, $d = 0.5$, timestepping using symplectic and symmetric implicit Runge-Kutta method \texttt{irk2} with step size 0.01.} 
	\label{fig:unstablekak}
\end{figure}

We can similarly construct a kink-antikink from two onsite kink structures (\cref{fig:unstablekak}, top left). The spectrum corresponding to this solution is shown in the top right panel of \cref{fig:unstablekak}, in which we see the split Goldstone modes on the real axis, confirming the instability of
this bound state. The corresponding Goldstone eigenfunctions are shown in the bottom left panel of \cref{fig:unstablekak}. Finally, we perform timestepping experiments for perturbations of this kink-antikink. If we perturb the two central nodes (both of which have values $u_n$ = 0) by a small amount, the solution develops oscillatory behavior about the neutrally stable intersite kink-antikink (\cref{fig:unstablekak}, bottom). The solution departs from the corresponding energy
maximum and performs oscillations around the nearby
energy minimum, namely the stable intersite kink-antikink
state. The energy $H(u(n, t))$ is again very well conserved as $t$ evolves (\cref{fig:unstablekak}), with a relative deviation from its initial value of less than $10^{-15}$. We note that for timestepping simulations on long time intervals, energy is very well conserved until boundary effects come into play. For $N=200$ and $N=400$ grid points, these boundary effects occur at approximately $t=200$ and $t=450$, respectively; enlarging the spatial grid delays these boundary effects.

\section{Conclusions and future challenges}\label{sec:conclusions}

In this paper, we used Lin's method to construct multi-kink solutions to the discrete Klein-Gordon equation by splicing together an alternating sequence of kink and antikink solutions with small amplitude remainders. These solutions exist as long as the distances between adjacent kinks and antikinks are sufficiently large. We then used Lin's method again to reduce the eigenvalue problem for multi-kinks to an effective, low-dimensional matrix 
eigenvalue equation. This matrix equation is different for each eigenvalue of the primary kink, and we find that for an $m$-component multi-kink, there are $m$ eigenvalues near each eigenvalue of the primary kink. Most notably, if the spectrum of the primary kink is imaginary, the spectrum of a multi-kink constructed from these primary kinks can be shown through this
explicit calculation to be imaginary as well. These eigenvalues can be computed numerically, and the result is in good agreement with the theory; this approach explicitly illustrates the 
spectral stability of these multi-soliton solutions.

Further avenues of research include exploring multi-kinks in other models. One such model is the Ablowitz-Ladik type discretization of the $\phi^4$ model (AL-$\phi^4$) \cite{kevrekidis2003}
\begin{equation}
	\ddot{u}_n = \frac{1}{h^2}(\Delta_2 u)_n + 2 u_n - u_n^2 (u_{n+1}+u_{n-1}),
\end{equation}
which has an exact static kink solution $u_n = \tanh(a n + \xi)$, where $a = \frac{1}{2} \cosh^{-1}[(1+h^2)/(1-h^2) ]$ and $\xi$ is arbitrary. Another is an alternative discretization of the sine-Gordon equation \cite{Barashenkov2008}
\begin{equation}
	\ddot{u}_n \cos \left( \frac{u_{n+1} - u_{n-1}}{4} \right) 
	= \frac{4}{h^2} \sin \left( \frac{u_{n+1} - 2 u_n + u_{n-1}}{4} \right)
	- \sin \left( \frac{u_{n+1} + 2 u_n + u_{n-1}}{4} \right),
\end{equation}
which has an exact static kink solution
\begin{align}
u_n(t) = 4\arctan\left[ \exp \left(k n - \frac{v t}{\sqrt{1 - v^2}} \right) \right],
\end{align}
where $k$ is defined implicitly by
\begin{align}
\sinh \left(\frac{k}{2}\right) &= \frac{1}{\sqrt{1 - v^2}}\frac{h}{\sqrt{4 - h^2}} && -1 < v < 1.
\end{align}
Here, $h$ is the lattice spacing, which is connected to our parameter
$d$ via $h=1/\sqrt{d}$, while $v$ denotes the kink's speed.
Although we do not expect these equations to have kink-antikink equilibrium solutions, since the stable and unstable manifolds of the two equilibria do not intersect transversely, we expect that a kink-antikink state will be an {\it approximate} equilibrium solution. We may then be able to explain the time evolution of the kink-antikink dynamics in terms of interactions between their exponentially decaying tails. 
In this context, it would be interesting to explore the
potential mathematical relevance of the corresponding
spectral calculation.
While we have already mentioned
that similar findings to the ones presented herein
are applicable to other Klein-Gordon models, such as,
e.g., the discrete $\phi^4$ model, it would be of
interest to study such ideas in more complex model variants 
such as the sine lattice (SL) studied in the work
of~\cite{Takeno1986}. 
Generalizing such ideas to the context of higher-dimensional
Klein-Gordon models where the kinks are typically also 
dynamically robust would be another direction of interest.
{We could also explore the time evolution using initial conditions which have the same asymptotics at $\pm \infty$ as a multi-kink.}
Finally, it remains open to explore whether similar techniques 
can be used to study multi-site breathers in discrete 
Klein-Gordon lattices.

\appendix

\section{Proof of \texorpdfstring{\cref{th:KaKexists}}{Theorem 1}}\label{sec:proof1}

The proof is an adaptation of the proofs of Theorems 1 and 3 in \cite{Parker2020}. First, we rewrite the system as a fixed point problem. Expanding $F(u)$ in a Taylor series about $c_i K(n)$, we get
\begin{align}\label{eq:Feq0}
F(U_i^\pm(n)) &= F(c_i K(n) + \tilde{U}_i^-(n)) = 
D F(c_i K(n)) \tilde{U}_i^\pm(n) + G(\tilde{U}_i^\pm(n)),
\end{align}
where $G(\tilde{U}_i^\pm(n)) = \mathcal{O}(|\tilde{U}_i^\pm|^2)$ with $G(0) = 0$ and $DG(0) = 0$. Since $f(u)$ is odd, $f'(u)$ is even, thus $D F(c_i K(n)) = D F(K(n))$, and equation \cref{eq:Feq0} becomes
\begin{align}\label{eq:Feq1}
F(U_i^\pm(n)) &= 
D F(K(n)) \tilde{U}_i^\pm(n) + G(\tilde{U}_i^\pm(n).
\end{align}
Since the pieces $U_i^\pm$ in the ansatz \cref{eq:Upiecewise} must match at their endpoints, we obtain the following system of equations for the remainder functions $\tilde{U}_i^\pm$
\begin{align}
\tilde{U}_i^\pm(n+1) &= D F(K(n)) \tilde{U}_i^\pm(n) + G(\tilde{U}_i^\pm(n)) \label{eq:Wsystem1} \\
\tilde{U}_i^+(N_i^+) - \tilde{U}_{i+1}^-(-N_i^-) &= c_{i+1} K(-N_i^-) - c_i K(N_i^+) \label{eq:Wsystem2} \\
\tilde{U}_i^+(0) - \tilde{U}_i^-(0) &= 0.\label{eq:Wsystem3}
\end{align}
Let $\Phi(m, n)$ be the evolution operators for the linear difference equation 
\[
V(n+1) = D F(K(n)) V(n).
\]
By the stable manifold theorem, $|K(n) - S^+| \leq C r^{-|n|}$ for $n \geq 0$ and $|K(n) - S^-| \leq C r^{-|n|}$ for $n \leq 0$, thus $| DF(K(n)) - DF(S^\pm)| \leq C r^{-|n|}$. Since $S^\pm$ are hyperbolic fixed points, we can decompose the evolution operator $\Phi(m, n)$ in exponential dichotomies on $\Z^\pm$ by \cite{Parker2020}*{Lemma 2}. The proof then follows that of \cite{Parker2020}*{Theorems 1 and 3}. Briefly, we write equation \cref{eq:Wsystem1} in fixed-point form using the discrete variation of constants formula together with projections on the stable and unstable subspaces of the exponential dichotomy. As long as $N$ is sufficiently large, we use the implicit function theorem to solve for the remainder functions $\tilde{U}_i^\pm$ as well as the matching conditions \cref{eq:Wsystem2} and \cref{eq:Wsystem3}. Since $W^u(S^-)$ and $W^s(S^+)$ intersect transversely, we have the decomposition $\R^2 = T_{K(0)}W^u(S^-)\oplus T_{K(0)}W^s(S^+)$, thus, as in the proof of \cite{Parker2020}*{Theorem 3}, solving \cref{eq:Wsystem3} does not involve jump conditions. The estimates \cref{eq:Uestimates} follow from the proof of \cite{Parker2020}*{Theorem 3} (see in particular \cite{Parker2020}*{Lemma 4}).

\section{Proof of \texorpdfstring{\cref{th:stability}}{Theorem 2}}\label{sec:proof2}

For the multi-kink solution $U(n)$, we will take an ansatz which is a piecewise perturbation of $V_0(n)$. (This ansatz is suggested in \cite{Sandstede1998}*{Section 7}). Let
\begin{equation}\label{eq:Viansatz}
V_i^\pm(n) = s_i c_i V_0(n) + W_i^\pm(n),
\end{equation}
where $s_i \in \R$, $s = (s_1, \dots, s_n)$, and $c_i = (-1)^{i+1}$. Substituting this into \cref{eq:EVPdyneq} and simplifying using \cref{eq:V0eq}, the eigenvalue problem becomes
\begin{equation}\label{eq:Weq1}
\begin{aligned}
W_i^\pm(n+1)
&= DF(K(n)) W_i^\pm(n) + \omega_0 B W_i^\pm(n) \\
&\qquad + [G_i^\pm(n) + (\omega - \omega_0) B](s_i c_i V_0(n) W_i^\pm(n)),
\end{aligned}
\end{equation}
where
\begin{equation}\label{eq:Gipm}
G_i^\pm(n) = DF(U_i^\pm(n)) - DF(K(n)),
\end{equation}
and we used the fact that $DF(c_i K(n)) = DF(K(n))$ from the previous section. Using \cref{eq:Aomegaeq}, we rewrite \cref{eq:Weq1} as
\begin{align}\label{eq:Weq2}
W_i^\pm(n+1)
&= A(n; \omega_0) W_i^\pm(n) + [G_i^\pm(n) + (\omega - \omega_0) B](s_i c_i V_0(n) + W_i^\pm(n)).
\end{align}

In addition to solving \cref{eq:Weq2}, the eigenfunction must satisfy matching conditions at $n = \pm N_i$ and $n = 0$. Thus the system of equations we need to solve is
\begin{equation}\label{eq:eigWsystem1}
\begin{aligned}
& W_i^\pm(n+1)
= A(n; \omega_0) W_i^\pm(n) + [G_i^\pm(n) + (\omega - \omega_0) B](s_i c_i V_0(n) + W_i^\pm(n))\\
& W_i^+(N_i^+) - W_{i+1}^-(-N_i^-) = S_i s \\
& W_i^+(0) - W_i^-(0) = 0,
\end{aligned}
\end{equation}
where
\begin{equation}\label{defDid}
S_i s = s_{i+1} c_{i+1} V_0(-N_i^-) - s_i c_i V_0(N_i^+).
\end{equation}

Let $\Phi(m, n)$ be the evolution operator for the variational equation 
\begin{equation}\label{eq:vareq1}
	V(n+1) = A(n; \omega_0) V(n).
\end{equation}
Since we have the exponential decay \cref{eq:A0decay}, by \cite{Parker2020}*{Lemma 2}, we can decompose the evolution operator $\Phi(m, n)$ in exponential dichotomies on $\Z^\pm$. In particular, we have the estimates 
\begin{equation}\label{eq:dichotomyest}
\begin{aligned}
&|\Phi_+^s(m, n)| \leq C r_0^{-(m - n)} && 0 \leq n \leq m \\
&|\Phi_+^u(m, n)| \leq C r_0^{-(n - m)} && 0 \leq m \leq n \\
&|\Phi_-^s(m, n)| \leq C r_0^{-(m - n)} && n \leq m \leq 0 \\
&|\Phi_-^u(m, n)| \leq C r_0^{-(n - m)} && m \leq n \leq 0 \:,
\end{aligned}
\end{equation}
for the evolution operator on the stable and unstable subspaces of the exponential dichotomy. Furthermore, letting $E^{s/u}$ be the stable and unstable eigenspaces of $A_0$, and $P_0^{s/u}$ the corresponding eigenprojections, we have the decay rates
\begin{equation}\label{eq:eigendecay}
	|P_\pm^{s/u}(n) - P_0^{s/u}| \leq C r^{-|n|},
\end{equation}
where $P_\pm^{s/u}(n) = \Phi_\pm^{s/u}(n, n)$ are the stable and unstable projections for the exponential dichotomy. We note that the decay rates in \cref{eq:dichotomyest} involve the eigenvalues for the matrix $A_0$, whereas those from \cref{eq:eigendecay} come from \cref{eq:A0decay}. 

The variational equation \cref{eq:vareq1} has a bounded solution $V_0(n) = (v_0(n), v_0(n-1))$, and the adjoint variational equation 
\begin{equation}\label{eq:adjvareq1}
	Z(n) = A(n; \omega_0)^* Z(n+1)
\end{equation}
has a bounded solution 
\begin{equation}\label{eq:Z0}
	Z_0(n) = (-v_0(n-1), v_0(n)),
\end{equation}
which is orthogonal to $Z_0(n)$. By the stable manifold theorem,
\begin{equation}\label{eq:VZ0estimates}
	\begin{aligned}
		V_0(n) &\leq C r_0^{-|n|} \\
		Z_0(n) &\leq C r_0^{-|n|}.
	\end{aligned}
\end{equation}

As in \cites{Parker2020,Sandstede1998}, we will not in general be able to solve the system \cref{eq:eigWsystem1}. Instead, since $\R^2 = \spn\{ V_0(0), Z_0(0) \}$, we relax the third equation in \cref{eq:eigWsystem1} to obtain the system of equations
\begin{equation}\label{eq:eigWsystem2}
	\begin{aligned}
	& W_i^\pm(n+1)
	= A(n; \omega_0) W_i^\pm(n) + [G_i^\pm(n) + (\omega - \omega_0) B](s_i c_i V_0(n) + W_i^\pm(n))\\
	& W_i^+(N_i^+) - W_{i+1}^-(-N_i^-) = S_i s \\
	&W_i^+(0) - W_i^-(0) \in \C Z_0(n).
	\end{aligned}
\end{equation}
Using Lin's method, we will be able to find a unique solution to this system, but this solution will generically have $m$ jumps at $n = 0$ in the direction of the adjoint solution $Z_0(0)$. A solution to \cref{eq:eigWsystem2} is therefore an eigenfunction if and only if the $m$ jump conditions
\begin{equation}
	\xi_i = \langle Z_0(0), W_i^+(0) - W_i^-(0) \rangle = 0
\end{equation}
are satisfied. Since $G_i^\pm(n) = \mathcal{O}(\tilde{U}_i^\pm(n))$, we have the same estimates for $G_i^\pm(n)$ as we do for $\tilde{U}_i^\pm$ in \cref{eq:Uestimates}. The terms in \cref{eq:eigWsystem2} are the same as in \cite{Parker2020}, with the addition of a term of the form $G_i^\pm(n) V_0(n)$. To estimate that term, we combine the estimates \cref{eq:Uestimates} and \cref{eq:VZ0estimates}. For the interior pieces, we have
\begin{align*}
| G_i^-(n) V_0(n) | &\leq C r^{-N_{i-1}^-} r^{-(N_{i-1}^- + n)} r_0^n \\
&\leq C r^{-2 N_{i-1}^-} \left(\frac{r}{r_0}\right)^{-n} \\
&\leq C \left(\frac{r}{r_0}\right)^{\max\{N_1, \dots, N_m\}} r^{-2 N} \\
&\leq C r^{-2 N} ,
\end{align*}
where in the last step we used the fact that $r_0$ is fixed (it depends only on $\omega_0$), as are the distances $N_i$. The constant $C$ will be larger as $r_0$ decreases, which occurs as the eigenvalue $\lambda_0$ gets closer to the continuous spectrum boundary. The estimate for $G_i^+(n) V_0(n)$ is similar, and the estimates for the two exterior pieces are stronger. Thus we have the two uniform estimates
\begin{equation}\label{eq:unifest}
	\begin{aligned}
		\|G_i^\pm\| \leq C r^{-N} \\
		\|G_i^\pm V_0(n) \| \leq C r^{-2 N}.
	\end{aligned}
\end{equation}

The proof then follows as in \cites{Parker2020,Sandstede1998}. First, we write equation \cref{eq:Weq1} as a fixed point problem using the discrete variation of constants formula together with projections on the stable and unstable subspaces of the exponential dichotomy. Let $\delta > 0$ be small, and choose $N$ sufficiently large so that $r^{-N} < \delta$. Define the spaces
\begin{align*}
V_W &= \ell^\infty([-N_{i-1}, 0]) \oplus \ell^\infty([0, N_i])  \\
V_a &= \bigoplus_{i=0}^{n-1} E^u \oplus E^s \\
V_b &= \bigoplus_{i=0}^{n-1} \ran P_-^u(0) \oplus \ran P_+^s(0)\\
V_\lambda &= B_\delta(\omega_0) \subset \C \\
V_s &= \R^m.
\end{align*}
Then for
\begin{align*}
&W = (W_i^-, W_i^+) \in V_W  && a = (a_i^-, a_i^+) \in V_a \\
&\lambda \in V_\lambda  &&b = (b_i^-, b_i^+) \in V_b,
\end{align*}
the fixed point equations for the eigenvalue problem are
\begin{equation*}\label{fpeig}
\begin{aligned}
W_i^-(n) &= 
\Phi_s^-(n, -N_{i-1}^-) a_{i-1}^- + \sum_{j = -N_{i-1}^-}^{n-1} \Phi_s^-(n, j+1)
[G_i^-(j) + (\omega - \omega_0) B](s_i c_i V_0(j) + W_i^-(j))
 \\
&+ \Phi_u^-(n, 0) b_i^- - \sum_{j = n}^{-1} \Phi_u^-(n, j+1) 
[G_i^-(j) + (\omega - \omega_0) B](s_i c_i V_0(j) + W_i^-(j))\\
W_i^+(n) &= \Phi_s^+(n, 0; \theta_i) b_i^+ + \sum_{j = 0}^{n-1} \Phi_s^+(n, j+1) 
[G_i^+(j) + (\omega - \omega_0) B](s_i c_i V_0(j) + W_i^+(j))\\
&+ \Phi_u^+(n, N_i^+) a_i^+ - \sum_{j = n}^{N_i^+-1} \Phi_u^+(n, j+1) 
[G_i^+(j) + (\omega - \omega_0) B](s_i c_i V_0(j) + W_i^+(j)),
\end{aligned}
\end{equation*}
where $a_0^- = a_m^+ = 0$, and the sums are defined to be $0$ if the upper index is smaller than the lower index. We now follow the procedure in the proof of \cite{Parker2020}*{Theorem 2} and invert the system \cref{eq:eigWsystem2} in the following steps. First, we solve for the remainder functions $W_i^\pm(n)$ in the fixed point equations. Then, we use the second equation in \cref{eq:eigWsystem2} to solve for the  initial conditions $a_i^\pm$. Finally, we use the third equation in \cref{eq:eigWsystem2} to solve for the $b_i^\pm$. The procedure is identical to that in the proof of \cite{Parker2020}*{Theorem 2}, with both $\lambda$ and $\lambda^2$ replaced by $\omega - \omega_0$, $\tilde{H}_i^\pm(n)$ replaced by $c_i V_0(n)$, and the presence of additional terms of the form $G_i^\pm(n) V_0(n)$, which will be small using the second estimate in \cref{eq:unifest} and will therefore be incorporated into the remainder term. Thus, we obtain a unique solution to \cref{eq:Wsystem3} which generically has $n$ jumps in the direction of $Z_0(0)$. These jumps are given by 
\begin{equation}\label{eq:xieq}
\begin{aligned}
\xi_i &= \langle Z_0(N_i^+), P_0^u S_i s \rangle 
+ \langle Z_0(-N_{i-1}^-), P_0^s S_{i-1} s \rangle - (\omega - \omega_0) s_i c_i M + R(\omega - \omega_0)_i(s),
\end{aligned}
\end{equation}
where $M$ is the Melnikov sum
\[
M = \sum_{j = -\infty}^{\infty} \langle Z_0(j+1), B V_0(j)\rangle,
\]
and the remainder term has uniform bound
\begin{equation}\label{eq:Rbound}
R(\omega - \omega_0)_i(s) \leq
C \left( r^{-N} r_0^{-2N} + |\omega - \omega_0|)(r_0^{-N} + |\omega - \omega_0| )\right).
\end{equation}
For the terms in the Melnikov sum,
\[
\langle Z_0(j+1), B V_0(j)\rangle = \frac{1}{d} \left\langle (-v_0(j), v_0(j-1))^\mathrm{T} , (v_0(j), 0)^\mathrm{T} \right\rangle
= -\frac{1}{d}v_0(j)^2,
\]
thus it follows that
\begin{equation}\label{eq:M}
M = \sum_{j = -\infty}^{\infty} v_0(j)^2 = \| v_0 \|_{\ell^2},
\end{equation}
which is always positive. It follows from 
\cref{eq:eigendecay} that
\begin{align*}
P_0^u S_i s &= s_{i+1} c_{i+1} V_0(-N_i^-) + \mathcal{O}\left( r^{-N}r_0^{-N}\right) \\
P_0^s S_i s &= -s_i c_i V_0(N_i^+) + \mathcal{O}\left( r^{-N}r_0^{-N}\right).
\end{align*}
Substituting these into \cref{eq:xieq}, we have
\begin{equation*}
\begin{aligned}
\xi_i = \langle Z_0(N_i^+), &V_0(-N_i^-) \rangle s_{i+1} c_{i+1}
- \langle Z_0(-N_{i-1}^-), V_0(N_{i-1}^+) \rangle s_{i-1} c_{i-1} \\
&+ \dfrac{1}{d} (\omega - \omega_0) s_i c_i M+ R(\omega - \omega_0)_i(d),
\end{aligned}
\end{equation*}
where the bound on the remainder term is unchanged. Since $c_i$ and $c_{i+1}$ have opposite signs, this simplifies to 
\begin{align*}
\xi_i = \langle Z_0(N_i^+), V_0(-N_i^-) \rangle s_{i+1}
- \langle Z_0(-N_{i-1}^-), V_0(N_{i-1}^+) \rangle s_{i-1}
- \dfrac{1}{d} (\omega - \omega_0) s_i M
+ R(\omega - \omega_0)_i(s).
\end{align*}
Let
\begin{align*}
a_i &= \langle Z_0(N_i^+), V_0(-N_i^-) \rangle 
= v_0(N_i^+)v_0(-N_i^- - 1) - v_0(N_i^+ - 1)v_0(-N_i^-).
\end{align*}
Using equation \cref{eq:Z0},
\[
\langle Z_0(-N_i^-), V_0(N_i^+) \rangle 
= v_0(N_i^+ - 1)v_0(-N_i^-) - v_0(N_i^+)v_0(-N_i^- - 1) = -a_i,
\]
which we substitute above to obtain
\begin{align}\label{eq:xifinal}
\xi_i = a_i s_{i+1}
+ a_{i-1} s_{i-1}
- \dfrac{1}{d} (\omega - \omega_0) s_i M 
+ R(\omega - \omega_0)_i(s).
\end{align}
The $a_i$ can be simplified further based on symmetries of the eigenfunction $v_0(n)$. For even, intersite eigenfunctions, $v_0(-n) = v_0(n-1)$, thus we have
\begin{align}\label{eq:aievenintersite}
a_i = v_0(N_i^+)v_0(N_i^-) - v_0(N_i^+ - 1)v_0(N_i^- -1).
\end{align}
For odd, intersite eigenfunctions, $v_0(-n) = -v_0(n-1)$, and so
\begin{align}\label{eq:aioddintersite}
	a_i = v_0(N_i^+ - 1)v_0(N_i^- -1) - v_0(N_i^+)v_0(N_i^-).
\end{align}
This can be written as the matrix equation
\begin{equation}\label{eq:matrixAeq}
	\left( A - \frac{1}{d}(\omega - \omega_0)MI + R(\omega-\omega_0)\right)s = 0,
\end{equation}
where $A$ is defined by \cref{eq:matrixA}, which has a nontrivial solution if and only if
\begin{equation}\label{eq:matrixdeteq}
	E(\omega - \omega_0) = 
	\det \left( A - \frac{1}{d}(\omega - \omega_0)MI + R(\omega-\omega_0)\right)s = 0.
\end{equation}
Let $\{ \mu_1, \dots, \mu_{m} \}$ be the eigenvalues of $A$, which are real since $A$ is symmetric. They are also distinct by \cite{Jirari1995}*{Corollary 2.2.7}, since the eigenvalue problem $(A - \mu I) v = 0$ is equivalent to the Sturm-Liouville difference equation with Dirichlet boundary conditions
\begin{equation*}
\begin{aligned}
\nabla( p_j \Delta s_j ) + q_j &= \mu s_j && \qquad j = 1, \dots, m \\
s_0 &= 0 \\
s_{m+1} &= 0,
\end{aligned}
\end{equation*}
where $p_j = a_j$, $q_j = a_j + a_{j-1}$, $\Delta$ is the forward difference operator $(\Delta f)_j = f_{j+1} - f_j$ and $\nabla$ is the backward difference operator $(\nabla f)_j = f_j - f_{j-1}$. 

Following the proof of \cite{Parker2020}*{Theorem 5}, we can use the implicit function theorem to solve for $\omega - \omega_0$ to get the $m$ solutions
\begin{align*}
	\omega_j &= \omega_0 + \frac{d \mu_j}{M} + \mathcal{O}(r_0^{-3N}) && j = 1, \dots, m.
\end{align*}
Since $\mu_j = \mathcal{O}(r_0^{-2N})$, $\omega_0 < 0$, and $\omega$ must be real, $\omega_j < 0$ for sufficiently small $N$. The corresponding eigenvalues are $\lambda_j = \pm \sqrt{\omega_j}$, which are on the imaginary axis.\\

\paragraph{Acknowledgments}

This material is based upon work supported by the U.S. National Science Foundation under the RTG grant DMS-1840260 (R.P. and A.A.)
and DMS-1809074 (P.G.K.).

\bibliographystyle{amsplain}
\bibliography{DKG.bib}

\end{document}